\documentclass{amsart}
\usepackage{graphicx}
\usepackage{amscd}

\newtheorem{thm}{Theorem}[section]
\newtheorem{cor}[thm]{Corollary}
\newtheorem{lem}[thm]{Lemma}
\newtheorem{prop}[thm]{Proposition}

\newtheorem{eg}[thm]{Example}
\theoremstyle{definition}
\newtheorem{defn}[thm]{Definition}

\numberwithin{equation}{section}

 \newcommand{\starspace}{
\begin{center}
* $\qquad$ * $\qquad$ *
\end{center}
}

\newcommand{\Pp}{\mathcal P}
\newcommand{\Pkn}{{\mathcal P}^k (n)}

\newcommand{\nbhd}{neighborhood }

 \newcommand{\bP}{\mathbb P} 
\newcommand{\R}{\mathbb R}   
\newcommand{\C}{\mathbb C}

\begin{document}

\title[Monster towers]{Curve Singularities and   Monster / Semple Towers}
\author{A. L. Castro$^{\star}$, R. Montgomery$^{\ddagger}$}%
\address{Mathematics dept. at UCSC, 1156 High st. ,  Santa Cruz - CA 95064}%
\email{alcastro@ucsc.edu$^{\star}$\footnote{Corresponding author.},rmont@count.ucsc.edu$^{\ddagger}$}%
\thanks{A.C. was supported by a UCSC Graduate Research Mentorship Award in the academic year of 08-09.}%

\date{21/11/09}%


\maketitle
\tableofcontents
\addtocontents{toc}{\protect\setcounter{tocdepth}{1}}

\begin{abstract}
In earlier work, we introduced the `Monster tower', a tower of   fibrations  associated to planar curves.
      We constructed an algorithm for classifying  its points  with respect to the equivalence relation
generated by     the action of the contact pseudogroup on the tower.
 Here,   we construct the analogous tower for curves in  $n$-space.
(This tower is known as the Semple Bundle in Algebraic Geometry.)
The pseudo-group of diffeomorphisms of $n$-space acts on each level of the extended tower.
We take initial steps toward classifying points of this extended Monster tower under this pseudogroup action.
 Arnol'd's list of  stable simple curve singularities plays a central role in these initial steps.  We end with a list of open problems.
\end{abstract}
\maketitle
\tableofcontents
\section{Introduction and History.}

Earlier works \cite{monz:flag}, \cite{monz:monster}, \cite{LJ:semple}   constructed  a `Monster  tower' :
\begin{equation*}
 \ldots \to M_k \to M_{k-1} \to \ldots \to  M_1 \to M_0.
 \end{equation*}
 of manifolds $M_k$ associated to curves in the  plane $M_0$.
The maps $M_k \to  M_{k-1}$ are fibrations with fiber the projective line.
 Each $M_k$ is endowed with a rank $2$ distribution $\Delta_k$.
 Here, and throughout, ``distribution'' means sub-bundle of the tangent bundle.
 This tower is called the Semple Tower in algebraic geometry (see   \cite{LJ:semple}).
 The base   $M_0$ can be taken to be any analytic surface in place of the plane.
 The tower is constructed
  by    iterating  Cartan's prolongation procedure for distributions as described in \cite{bry:cartan}.
  By a `symmetry' at level $k$ we mean a local analytic diffeomorphism of $M_k$
  which maps $\Delta_k$ to itself.  The distribution $\Delta_1$
  is a contact distribution and its symmetries are the contact (pseudo) group.  A theorem of Backl\"und-Yamaguchi ( \cite{yama:cont} ) asserts that all symmetries at level $k > 1$
  arise  via prolongation  from the symmetries at  level 1.
   The central
  problem  addressed  in \cite{monz:monster}  was to classify  the orbits of this symmetry group at any level.
  We largely solved this problem by constructing    an algorithm for   converting  it to the well-studied  problem of
  classifying   finite jets of plane curve singularities, and using tools such as the
  Puiseux characteristic, well-known in that case.

In the present paper we    take the first steps towards
 generalizing this work from the plane   to n-space, $\C^n$.
 We   construct an analogous tower:
 \begin{equation*}
 \ldots \to \Pp^k (n) \to \Pp^{k-1}(n) \to \ldots \Pp^1 (n) \to \Pp^0 (n) = \C^n.
 \end{equation*}
for   curves in  $n$-space, $\C^n$.    The fibers are now   projective  spaces of dimension $n-1$.
  When $n =2$, we have $\Pp ^k (2) =M_k$ of above.
When $n > 2$, Backl\"und-Yamaguchi's theorem now
asserts that all  symmetries at level $k$ arise from  level $0$ where the symmetry (pseudo-) group
is $Diff(n)$, the pseudo-group of locally defined analytic diffeomorphisms of $\C^n$.
 We   solve  some first   occurring   instances of
 the  corresponding   classification problem by describing
      ``first occurring'' orbits that do not arise in  the planar case $n=2$
and by   classifying  the codimension 1 and codimension 2 singularities at any level, for any $n$.
{\bf We conjecture that the  problem of classifying simple stable (with increasing $n$) singularities of the Monster
 tower is  equivalent  to Arnold's classification \cite{arn:sing} of simple stable
 curve singularities.}  We verify
   the first instances of such a
 correspondence in the   course of classifying  codimension 1 and 2 singularities for
 the Monster tower.   We also get lower bounds on the number of
 orbits within $\Pp^k (3)$, for $k =2,3,4,5$, indicating many more orbits
 than in the planar case and discuss relations with the classical Enriques formula
 as obtained by \cite{LJ:semple}.

The  monster  tower $\Pkn$ is known as
the ``Semple Tower'' in algebraic geometry. It was introduced
by Semple (\cite{semp:tower}).   See in particular the discussion in  (\cite{LJ:semple, ken:cnct}). The algebraic geometers typically take the  base  $\Pp^0 (n)$
of their tower to be n-dimensional projective space, or a general smooth n-dimensional variety, instead
of  our $\C^n$.

 The tower $\Pkn$  is  the universal embedding space for the  Nash Blow-ups of curves in $\C^n$.
  Nash blow up is an alternative to the usual blow up of algebraic geometry,
  in which the secants lines of  the usual blow-up are replaced by tangent lines.
   (See pp. 412, 3rd paragraph of
op. cit. \cite{spi:nash} for a history and Nash's original manuscript \cite{nash:bu}. ( See also \cite{har:first}, esp. pp. 219-221. )
 In this paper, we lean towards  use of the
 Cartan language so refer to the   $k$-th Nash blow-up of a
 curve as its ``$k$-th prolongation''.  The $k$-th prolongation of a curve in $\C^n$  is an
 integral curve for $\Delta_k$ in $\Pkn$.
 A   theorem of Nobile \cite{nob:nash}
  asserts that for  sufficiently large $k$, the  $k$-th prolongation  of a singular algebraic or analytic  curve is smooth.\\

\paragraph{{\em Note.}} We would like to point out that P. Mormul has several papers (\cite{morm:moduli,morm:arxiv}) of a kind of parallel nature
concerning the classification of points in $\Pkn$, $n > 2$, interpreted as flags of special Goursat distributions.

\starspace

\paragraph{ {\bf Acknowledgements.} } A.C. would like to thank J. Castellanos and E. Casas-Alvero (Madrid and Barcelona resp., Spain) for  teaching us about the subtleties of the various discrete invariants of space curves , and  P. Mormul (Warsaw, Polland) for clarifying some questions regarding appearance of moduli in the classification problem of (special) Goursat flags. Thanks  to G. Kennedy (Ohio, US) for providing us with some reprints of his work , and from where we could learn more about the Semple tower. Special acknowledgments go to Misha Zhitomirskii (Haifa, Israel), our singularity theory mentor, and without whom this project would have never taken off.

\section{Preliminaries.  Construction}

All manifolds, curves,  maps, etc. are assumed   analytic.
We will work over $\C$ instead of $\R$ but all results hold
for the real case also.   As in \cite{LJ:semple},  one could
also work over any algebraically closed field of characteristic zero in place of $\C$.

\subsection{Prolongations}
Let  $D$ be a rank $m$  distribution on a manifold $Z$.
Viewing $D$ as a vector bundle over $Z$, we   form its
projectivization $\pi: \bP D \to Z$,  the fiber bundle whose
fiber over $z$  is the $m-1$ dimensional projective space $\bP (D(z))$.
A point of $\bP D$ is then a pair $(z, \ell)$ with $z \in Z$
and $\ell \subset D(z)$ a line in the $m$-dimensional vector
space $D(z) \subset T_z Z$.  Define
a distribution $D^1$ on $\bP D$ by
setting $D^1 (z, \ell ) = d \pi_{(z, \ell)} ^{-1} (\ell)$.
Since  $\pi$ is a submersion  with $(m -1)$-dimensional fibers,
and since   $\ell$ is one-dimensional
we have $rank(D) = rank(D^1) = m$.

Let  $\gamma \subset Z$ be a non-constant parameterized  integral curve for $D$.
Its prolongation, denoted $\gamma^1$,
is  the  integral curve for $D^1$ defined by
$\gamma^1 (t) = (\gamma(t),span ( d \gamma /dt))$ at regular points $t$
for a local parameterization of $\gamma$.  If $t=t_0$ is not regular,
we define $\gamma^1 (t_0)$ by taking the limit $lim_{t \to t_0} \gamma^1 (t)$ of
regular points $t \to t_0$.  Our analyticity assumption implies
that this limit is well-defined and that the resulting  curve $\gamma^1 (t)$ is
analytic everywhere.   For a proof, see \cite{monz:monster}, pp. 14.

Let $\phi$  be a symmetry of $(Z,D)$, meaning a diffeomorphism of $Z$
for which $\phi_* D = D$. Its prolongation,
  $\phi^1$  is the symmetry of $(\bP D, D^1)$
defined by $\phi^1 (z, \ell) = (\phi(z), d \phi_z (\ell))$.
We have   $(\phi \circ \gamma)^1 = \phi^1 \circ \gamma^1$
for $\gamma$ an integral curve for $D$.

This prolongation construction is due to \'E. Cartan.  We learned of it from R. Bryant (\cite{bry:cartan} ) .

\starspace

\paragraph{{\it Note.} } We have   described   the ``rank 1''  prolongation of $D$,
``rank 1'' in that this prolongation  is associated to lines, and so integral curves.   For each $r < m=  rank(D)$ there is a ``rank $r$'' prolongation      associated to $r$-dimensional integral
submanifolds for $D$.  The  $r > 1$ prolongations are  significantly
more complicated than the $r=1$ prolongations due to the need
to   account for the equality of mixed partial derivatives in forming integral submanifolds.
 Look for  ``integral elements'' in \cite{bry:extdiff}  for
some details, and see the recent paper by Shibuya-Yamaguchi \cite{shib:drapeau}.

\subsection{Building the tower} \label{subsection:prol}

Start with complex  $n$-space  $\C^n$ endowed with
its tangent bundle $\Delta_0 = T \mathbb \C ^n$ as distribution.
Prolong  to get the manifold  $$\mathcal{P}^1 (n) = \bP \Delta_0 \cong \C^n \times \C \bP ^{n-1}$$ with distribution $\Delta_1$.
Repeat.  After $k$ iterations  we obtain a  manifold $\mathcal{P}^k(n)$
endowed with the rank $n$ distribution $\Delta_k = (\Delta_{k-1}) ^1$.  Topologically,
  $$\mathcal{P}^k (n)= \mathbb C^n \times \bP^{n-1} \times \ldots \to \bP^{n-1}$$ ($k$ copies
of   projective space), and  the projection $$\mathcal{P}^k(n) \to \mathcal{P}^{k-1}(n)$$   projects out
the last factor.

\begin{defn}  The {\em Monster  tower} for curves in $n$-space  is   the sequence
of manifolds with distributions, $(\mathcal{P}^k(n), \Delta_k)$,
together with the fibrations  $$\ldots \to \mathcal{P}^k(n) \to \mathcal{P}^{k-1}(n) \to \ldots
\to \mathcal{P}^0(n) = \C^n.$$
\end{defn}
We write $\pi_{k, i}: \Pkn \to \Pp ^{i} (n)$, $i \le k$
for the projections.

\medskip

 The group of analytic diffeomorphisms of $\mathbb C^n$ acts on the Monster tower
by prolongation, preserving the levels, as indexed by $k$, the  distributions
$\Delta_k$ and the fibrations.  To avoid
 restrictions arising  from  convergence and domain of
definition issues  and to allow more flexibility
in the equivalence relations we will work with pseudogroups and germs instead
of globally defined diffeomorphisms.   Let $Diff(n)$
be the pseudogroup of analytic diffeomorphisms of $\C^n$.   If $\psi \in Diff(n)$
has domain $U$ and range $V$ then its $k$-th prolongation
$\psi^k$  will have domain $\pi_{k,0}^{-1}(U)$
and  range $\pi_{k,0}^{-1}(U)$.

 \begin{defn} We say that two points $p, q \in \Pp^k (n)$ are equivalent, in symbols
$p \sim q$, if
there is a diffeomorphism germ $\psi \in Diff(n)$ such that
$\psi ^k (p) = q$.
\end{defn}

\medskip

{\bf Classification Problem.}  Classify the resulting equivalence classes.

\medskip

Conceptually, it  is often  simpler to fix the base points $p_0 = \pi_{k, 0} (p)$ and $q_0=
 \pi_{k, 0} (q)$
to be $0 \in \C^n$.
 Then we can replace the pseudogroup  $Diff(n)$
by the {\em honest} group $Diff_0(n)$ of germs of diffeomorphisms of $\C^n$ mapping $0$ to $0$.
The classification problem is  then replaced by the   problem
 of classifying the   orbits  for this action on the
fiber $\Pp^k (n)_0 = \pi_{k,0} ^{-1} (0) \subset  \Pp^k (n)$ over $0$.

\medskip

\section{Language. Results.}

\subsection{The Curve Approach.} Take a non-constant  curve  $\gamma(t)$ in   $\mathbb C^n$  with $\gamma(0) = 0$.
Prolong it repeatedly
to form the sequence of curves $\gamma^1, \gamma^2, \ldots$ with   $\gamma^k (t)$
a curve in $\mathcal{P}^k(n)$,
 integral for $\Delta_k$ and $\pi_{k, i} \circ \gamma^k = \gamma^i$.  Since
 $(\phi \circ \gamma)^k (t) = \phi^k \circ \gamma^k (t)$
for any $\phi \in Diff_0 (n)$ we see  that
\begin{equation} \label{CurveDict}
\text{ if }  \\   p = \gamma^k (0) \\
\text{  then }  \\ \phi^k (p) = (\phi \circ \gamma)^k (0).
\end{equation}
 This observation
suggests that we  approach our  classification problem by turning it into the well-studied
  classification problem for curve germs.    {\em This curve approach to the problem
for  the case   $n =2$  was very successful \cite{monz:monster}.}

\begin{defn} [preliminary]  For $p \in \mathcal{P}^k(n)$,
suppose that there is a non-constant curve germ $\gamma$ in $\C^n$
for  which $\gamma^k (0) = p$ and that $\gamma^k$  is immersed.
Then we say that $\gamma$ realizes $p$.
\end{defn}
Although we have required that $\gamma^k$ be immersed,
it  is essential to the definition that we allow  $\gamma$ to be singular.

Let   $\tau: (\C, 0) \to (\C, 0)$ be  a non-constant map.
Observe that   $$( \gamma \circ \tau )^k (t) = \gamma^k  \circ \tau $$
which in particular asserts that if  $\gamma$ realizes $p$,
so does any re-parameterization of $\gamma$.
In the standard equivalence relation for curves,  often called
``RL'' (for right-left) equivalence we insist that
$\tau$ is a local diffeomorphism (= reparameterization).
\begin{defn} \label{RL} Equivalence of curve germs.
We will say   two  curve germs $\gamma$,  $\sigma$ in $\C^n$ are equivalent,
in symbols $\gamma \sim \sigma$,  if and only if there is a diffeomorphism germ
 $\psi \in Diff(n)$ and a reparameterization  germ $\tau \in Diff_0 (1)$
 ( thus $\tau(0) = 0; d \tau/ dt (0) \ne 0$ )
of $(\C, 0)$   such that
$\sigma = \phi \circ \gamma \circ \tau ^{-1}$.
\end{defn}

Recall  that two points $p, q \in \Pp^k (n)$ are equivalent , in symbols
$p \sim q$, if
there is a diffeomorphism germ $\psi \in Diff(n)$ such that
$\psi ^k (p) = q$.
 We have just seen that if $p$ and $q$ are points at the same level,
realized by curves $\gamma, \sigma$ then  $\gamma  \sim \sigma$ implies
that  $p \sim q$.   {\bf This fact is the heart of the curve approach.}

From the singularity viewpoint, the simplest of all curve germs are
the immersed curves -- those with $d \gamma / dt  (0) \ne 0$,
 and they are all equivalent.  \begin{defn}  A  point $p$ of the Monster is called   a `Cartan point' if it can be realized by  an immersed  curve germ $\gamma$ in $\C^n$.
\end{defn}

The following theorem is well-known. It can be found in \cite{yama:cont} for example.
\begin{thm} \label{Cartan}
For $k > 1$ the Cartan points of $\mathcal{P}^k(n)$  form a dense open orbit  $C^k (n)$
whose complement is a hypersurface.
In a \nbhd of any Cartan point the distribution $\Delta_k$
is locally diffeomorphic to the canonical distribution on
the space $J^k (\C, \C^{n-1})$ of k-jets of analytic  maps from $\C$
to $\C^{n-1}$.
Every point of $\mathcal{P}^1(n)$ is a Cartan point.
\end{thm}
We reprove the theorem here, both for completeness,
and because in developing the proof we will develop needed tools.
Our proof is in two parts, one being the proof of Proposition \ref{CartanNV}
further on in this section, and
the other being Example 1 within
the subsection of section 4 on KR coordinates.

\begin{defn}  The  singular locus at level $k$ is the set
$\mathcal{P}^k(n) \setminus C^k (n)$. A singular point is a point of the
singular locus.
\end{defn}

Theorem \ref{Cartan}  suggests that we will need singular curves to realize  singular points of the Monster.  The first occuring singularity in any list of singular curves are the $A_k$ singularities.  For even $k$, these
are single-branched,
represented  as the parameterized curve
\begin{equation}
\label{A2k}
x_1 = t^2,  x_2 = t^{2k +1}    \hskip 2cm ( \text{ the } A_{2k} \text{ curve } )
\end{equation}
To include   curve ( \ref{A2k}) in $\mathbb C^n$, $n > 2$,   set the remaining coordinates to zero:
$x_3 = \ldots  =  x_n = 0$.  Any curve diffeomorphic to a reparameterization of  the $A_{2k}$ curve is called
`an $A_{2k}$ singularity'.
Computations  done in the next section (`Example: the $A_{2k}$ singularity',  soon after KR coordinates are introduced)
 show that for $j \ge k$, the j-fold prolongation of
an $A_{2k}$ singularity is immersed.

\begin{thm} \label{codim1} The points in $\Pp ^j (n)$   realized
by    the $A_{2k}$ singularity (eq. \ref{A2k})   for  $j \ge k > 1$ and any $n \ge 2$
   form a     smooth quasiprojective hypersurface lying in the singular locus $\Pp ^j (n) \setminus C^j (n)$.   The
  union of these $A_{2k}$ points  over this  range of $k$ is   open and
dense within the singular locus at level $j$.
\end{thm}

We prove this theorem in section 4.  To proceed to further results
we need some language from singularity theory.
\starspace
\paragraph{\bf Synopsis and terminology.}
 A {\em singularity  class} for  $\Pkn$ is  a subset of $\Pkn$ which  is  invariant under the
  $Diff(n)$ action.    If that singularity class forms a subvariety, then
  its  codimension is the usual codimension of this subvariety,
  as a subvariety of $\Pkn$.   A point, or a singularity class, is   called  {\em simple}  if  it  is contained in a \nbhd which  is the union of  a finite number of
$Diff(n)$ orbits.

The `order' or `multiplicity' of an analytic function germ  $f(t) = \Sigma a_i t^i$
is the smallest integer $i$ such that $a_i \ne 0$. We write $ord(f)$
 for this (non-negative) integer.
The {\em multiplicity} of an analytic  curve germ $\gamma: (\C, 0) \to (Z, p)$
 in a manifold $Z$  is the minimum of the orders of its coordinate functions $\gamma_i (t)$
 relative to any coordinate system vanishing at $p$.  This multiplicity, denoted
 $mult(\gamma)$ is independent of a parameterization of $\gamma$ and of
 choice of a vanishing coordinate system.   (Warning: Other uses of the term  `multiplicity'
 applied to curve singularities abound in the singularity theory literature.   We are following
 the use of the word  as found in Zariski \cite{zar:mod} or Wall
 \cite{wall:curv}.   More precisely, every single-branched plane curve singularity
 $f(x,y) = 0$ with $f(0,0) = 0$ can be well- parameterized as a curve $\gamma (t) = (x(t), y(t))$
 and what we have called the `multiplicity of $\gamma (t)$ at $t =0$ coincides
 with Zariski's multiplicity in this case. )

  The {\em embedding dimension} of a curve in the manifold $Z$ is
  the smallest integer $d$ such that $c$ lies in a smooth  $d$-dimensional  manifold $\Sigma \subset Z$.
  Thus, embedded curves have embedding dimension $1$, while $A_{2k}$ curves have embedding dimension
  2. A curve germ with embedding dimension at least 3 is called a  {\em spatial} curve.

  A curve germ $\gamma: (\C, 0) \to Z$  is called simple
if, for all  integers $N$ sufficiently large the  $N$-th jet $j^N \gamma$ of $\gamma$ is contained in
a \nbhd which is covered by a finite number of RL-equivalence classes.

We summarize the    classification results thus  far. Fix the level $k_0$.
The Cartan points at that level form a single orbit, which is open and dense and whose points are
realized by immersed (multiplicity $1$) curves.
The generic  singular   points
at that level are realized by the $A_{2k}$ singularities for $k \le k_0$.
For each $k$ these realized points form   a single orbit of  codimension 1.
The
realizing curves -- the $A_{2k}$ singularity -- has   multiplicity 2 and is strictly planar. In both cases, the   representative curves are  simple, as are the  points they realize.
These facts suggest the following refinements to  the classification problem.

\medskip

{\bf Refined Classification Questions.}
Do the  codimension 2 singularities  within the Monster correspond to multiplicity 3 curves?


Do simple singularities of the Monster always correspond to
simple curve singularities?

At what level, and  in what codimension, do non-planar singularities first occur?

  We answer the first question with theorem \ref{thm:cod2class}
  and the last question in Corollary (Cor. \ref{1stOccurring} ).
  Theorem \ref{thm:cod2class}  and the discussion following it,
  together with the final section, gives partial answers to the second question.

\medskip

\subsection{Verticality. Baby Monsters. RVT coding. }

As with   any smooth fiber bundle, we have the notion
of  the `vertical space' for the fibration $\mathcal{P}^k(n) \to \mathcal{P}^{k-1}(n)$.

\begin{defn}  The vertical space at $p$ is the  linear subspace
\begin{equation*}
V_k (p) = ker(d \pi_{k, k-1} (p)) \subset T_p \Pp^k (n).
\end{equation*}
A vector $v \in T_p \Pp^k (n)$
or line $\ell  \subset T_p \Pp ^k (n)$ is called vertical if $v \in  V_k (p)$
or $\ell \subset V_k (p) $.
\end{defn}
Since the vertical spaces are the  tangent spaces to the fibers of $\pi_{k, k-1}$
and since the symmetry group $Diff(n)$ maps fibers to fibers,
we have that symmetries send vertical spaces to vertical spaces.
Moreover, the fibers are integral manifolds of the distribution $\Delta_k$,
so that
$$V_k (p) \subset \Delta_k (p).$$
 \begin{defn} The point $p = (p_{k-1}, \ell)$ of the Monster at level $k$, $k > 1$  is called
`vertical' if the line $\ell$ is  a vertical line at  level $k-1$. Otherwise, that point is called `non-vertical'.
Every point at level 1 or 0 is considered to be non-vertical.
\end{defn}

The following proposition is basic to understanding the  Cartan points
and the overall structure of the Monster tower.
\begin{prop} \label{CartanNV} Let $p \in \Pkn$ and write $p_i =  \pi_{k, i} (p)$, $i \le k$. The point $p$ is Cartan if and only if
none of the $p_i$'s below $p$ is vertical.
\end{prop}

{\bf  Proof:}  Postponed for a few pages.

\medskip

Note that this proposition contains  the bulk of    theorem \ref{Cartan}. All that remains to prove of theorem \ref{Cartan} is the
assertions regarding equivalence between points of
the jet bundle for maps to $\C^{n-1}$, and that is done in the first example of the next section.

Prolongation can be applied to any manifold $F$  in place of $\mathbb C^n$.
We let $\Pp ^0 (F) = F$, with its tangent bundle $\Delta_0 ^F = TF$ as distribution.
The prolongation of $(F, \Delta_0 ^F)$ is $\Pp ^1 =   \bP TF$,
with its canonical rank $m =dim(F)$ distribution $\Delta_1 ^F = (\Delta_0 ^F)^1$,  etc.
Locally $(\Pp ^k (F), \Delta_k ^F)$ is  analytically diffeomorphic, as a manifold equipped with a distribution,
to $(\Pp ^k (m), \Delta_k)$.
If $(Z, D)$ is a manifold with distribution
and $F \subset Z$ is an integral submanifold of $D$ then
$\Pp ^1  (F) \subset \bP D$, and $\Delta_1 ^F = D^1 \cap T (\Pp ^1  (F))$
over   $\Pp ^1  (F)$.

We apply these considerations to the fiber $F_k (p) : = \pi_{k, k-1} ^{-1} (p_{k-1}) \subset \Pkn$ through the point $p$
at level $k$.
(As above, $p_{k-1} =  \pi_{k, k-1} (p)$.)  The fiber  is an $(n-1)$-dimensional  integral submanifold for $\Delta_k$.
Prolonging, we get $\Pp^1 (F_k (p)) \subset \Pp ^{k+1} (n)$,
together with its distribution,  $\delta_k ^1 = \Delta_1^{F_k (p)}$; that is,
$$ \delta_k ^1 (q)  =  \Delta_{k+1} (q) \cap T_q (\Pp^1 (F_k (p))   $$
 a hyperplane within $\Delta_{k+1} (q)$,  for $q \in \Pp^1 (F_k (p))$.
Iterating,  we obtain embedded   submanifolds
$$\Pp^j (F_k (p)) \subset \Pp ^{k+j} (n), $$
together with hyperplanes $\delta_k ^j (q) \subset \Delta_{k + j} (q)$
for $q \in \Pp^j (F_k (p))$.

\begin{defn}   We call the tower $(\Pp^j (F_k (p)) , \delta_k ^j)$, $j = 0, 1, \ldots$
the baby Monster through $p$.
\end{defn}

\begin{defn}
A {\em critical hyperplane} at level $k$ through a point $p$ is any one
of the hyperplanes   $\delta_i ^j (p)  \subset \Delta_k (p)$  (with $i + j = k$)
associated to the prolongation of a fiber through a point in the tower under $p$.

A direction $\ell \subset \Delta_k (p)$ is called
{\em regular} if it does not lie in any critical hyperplane.
A direction is called {\em critical} if it does lie in a critical hyperplane,
and is called {\em  `tangency' } if that critical hyperplane is not the vertical hyperplane.
 (A rationale for the tangency terminology can be found in   \cite{monz:monster}.)

A point is called regular, critical, vertical, or tangency, depending on
whether the corresponding line one level down is regular, critical, vertical,
or tangency.

An integral curve is called `regular' if it is tangent
to a regular direction.
\end{defn}

{\bf Warning:} a line can lie in more than one critical hyperplane.

\medskip

Note the vertical hyperplane $V_k (p)$
is itself  a critical hyperplane, being of the form $\delta_k ^0$.

\begin{thm} \label{regNonEmpty}Through every point there passes a regular integral direction.
\end{thm}

{\bf Proof:}  If $p$ is at level $k$ then there are at most $k-1$ critical planes
through $p$, one for each point in the tower  below $p$ besides $p_0$.
The complement of a finite collection of hyperplanes is open and dense.
Take $\ell$ to be any line in this complement. Q.E.D.

\subsubsection{RC codes and RVT codes.}  Every point $p$ is either regular or critical.
Let $p$ be a point and let $\{p_0, p_1, \ldots, p_k = p\}$ be the set of all projections of $p$ to lower levels. (Cf. notation of Proposition \ref{CartanNV}.)\\
The $RC$ code of $p$ is the word $w = w_1 \ldots w_k$ of length $k$ in
the letters $R$ and $C$ with $w_i = R$ if $p_i$ is regular and $w_i = C$ if
$p_i$ is critical.  The RC class of the word $w$ is the singularity class in $\Pkn$,
denoted $w \subset \Pkn$ by slight abuse of notation,  consisting of all those
points $p \in \Pkn$ having RC code $w$.

Note that $w_1 = R$ for any point, as all points at level 1 are regular.

\medskip

{\bf Example.}  According to Proposition \ref{CartanNV}, and the spelling rules,
 the Cartan points at level k
are those points whose RC code is $R^k = R \ldots R$ (k times).

\begin{prop} \label{RC} The codimension of the RC class $w$
is equal to the number of letters $w_i$ which are $C$'s.
An RC class $w$ adjoins an RC class $\tilde w$ -- meaning $w$ lies in
the closure of $\tilde w$ -- if and only if $\tilde w$ can be made into $w$
by replacing some of  the occurrences of the letter $C$ appearing in $w$ by the letter $R$.
\end{prop}

{\bf Proof.} The critical planes are hyperplanes within
the distribution and the condition that a line
lie in a given hyperplane is defined by a single equation.
To prove the second assertion, {\em realize that any critical point at any level
is the limit of regular lines passing through the same point one level down},
and hence $\tilde w$'s
closure contains $w$.
Q.E.D.

\medskip

{\bf Example.} If $c(t)$ is the  $A_{2k}$ curve,
then, we will compute at the end of section 4.2,  that  the RVT code of
the point  $c^{k +1 +s} (0)$  is   $R^k C R^s$,
corresponding to the fact that the orbit
of this curve has codimension 1. It adjoins the Cartan class
$R^{k+1 +s}$.

\medskip

\paragraph{{\bf RVT code.}} Change  our   alphabet by replacing $C$ by  $V$
or $T$.  Set $\omega_i (p) = V$ if $p_i$ is a vertical point, and set  $\omega_i (p) =T$
if $p_i$ is  a  a  critical point {\it which is not vertical}. In this way we associate  an RVT code
to each point, and an associated RVT class.
An RVT refinement $\omega$  of an RC code $w$
is any RVT code which becomes $w$  when all occurences of the letters V and T
are replaced by C.

\medskip

{\bf Example.}  If $w = RCCCRC$ then the  possible RVT refinements of $w$ are
$\omega = RVVVRV, RVVTRV,RVTVRV,RVTTRV$.

\medskip

{{\bf Spelling rules.}}   As discovered in the book \cite{monz:monster}, there are only two spelling rules for RVT words $w_1 \ldots w_k$.
The first rule is that  $w_1 = R$: every word starts with R.
The second rule is that the letter  T  cannot immediately follow the letter
 R, reflecting  the fact that for $p_i$ to be a  tangency point
it must lie in the baby monster of some point in the tower under $p_i$.

\medskip

\paragraph{{\bf Warning:  the class ``L''} } Critical planes are hyperplanes in an  $n$-dimensional space.
The intersection of two such planes will contain a line as soon as    $n >2$, and this line
represents a point one level higher.
For example, when  $n =3$  the intersection of the vertical
plane through $p$  and a critical plane    arising from a lower level
will contain a line.   One
is tempted to say  that the corresponding point, one level up, is
both a V and a T point.    We have chosen our terminology so
that it is labeled to be a V point, but a  special name is useful for such a point.
A  point  $p = (p_k , \ell)$
whose line  $\ell$ satisfies    $\ell \in V_k (p) \cap \delta^i _j (p)$
will be called an ``L'' point.
``L'' points   can be reached by
swinging a tangency line $\ell_t \in \delta^i _j$  around until it becomes vertical.
Thus  an L  point lies in   the intersection of  the closure of a  class ending in  T
with a class ending in a V.

Further refinements of the code are clearly possible, for example by indicating how many critical planes
a line lies in, at what level baby monster these planes originate, etc.  We leave these further
developments to interested parties.

\subsection{Points by Curves.}

We will associate to $p$ the collection $Germ(p)$ of all curve germs which
realize it regularly.

\begin{defn}  For $p$ a point at level $k$, write
$Germ(p)$ for  the collection of  all curve germs  $\gamma: (\C,0)  \to \C^n$ at $t=0$
whose $k$-th prolongation $\gamma^k$  is regular and passes through $p$:
$\gamma^k (0) = p$.  From now on, when we say ``$\gamma$ realizes $p$''
we mean that $\gamma \in Germ(p)$, so not only is $\gamma^k (0) = p$
and $d \gamma^k /dt |_{t = 0} \ne 0$, but also the span of $d \gamma^k /dt |_{t = 0} \ne 0$
is a regular direction.
\end{defn}

According to theorem \ref{regNonEmpty}, $Germ(p)$ is non-empty.
Since the  $k$-fold prolongation of the $k$-fold   projection of a
regular curve is the original  curve,
we have the following alternative description of $Germ(p)$:
$$Germ(p) = \{ \pi_{k,0} \circ \sigma :
\sigma  \text{ is a regular integral curve germ passing through }  p \}.$$

The following proposition is immediate:
\begin{prop} \label{1regularizing} Let $\sigma$ be a regular integral curve germ at level $k$.
\begin{itemize}
\item[(a)]
$\sigma^1 (t)$ is a regular integral curve germ at level $k+1$.
\item [(b)] If $\sigma(0)$ is a regular point, then the one-step projection
$\pi_{k, k-1} \circ \sigma$  of $\sigma$ is a regular integral curve germ.
\end{itemize}
\end{prop}

{\bf Proof of Proposition \ref{CartanNV} on Cartan points.}
We must prove a point is Cartan if and only if
its code is $R^k$.  Suppose $p$
is Cartan. Let $\gamma$ be an immersed curve representing
$p$.  Being immersed, all its prolongations $\gamma^1, \gamma^2, \ldots$
are immersed and regular, by (a) of the above proposition.  Thus all the points in the tower below $p$
are regular and so its code is $R^k$.   Conversely, suppose that all the
points in the tower below $p$ are regular.
Let $\sigma = \gamma^k$ be a regular curve passing through
$p$ and consider its one-step projection $\sigma_1 = \pi_{k , k-1} \circ \sigma$.
By (b) of the proposition, and the fact that $p_{k-1}$ is regular,
we have that $\sigma_1 = \gamma^{k-1}$ is regular.
Continuing to project we see that all the projections of $\sigma$
are regular. In particular $\gamma$ is immersed, and so $p$ is Cartan.
Q.E.D.

\medskip

Later on we  will need to know that curves in $Germ(p)$
are well-parameterized.

\begin{defn}\label{def:wellparameterized}
A curve is called well-parameterized if it
has a representative for which the map $t \to \gamma (t)$
is one-to-one.
\end{defn}

Equivalently,  a curve germ
 $\gamma$ is {\it not} well-parameterized if and only if
we can express $\gamma = \sigma \circ \tau$
for some other curve germ $\sigma: (\C, 0) \to \C^n$
and some {\it non-invertible}  germ $\tau: (\C, 0) \to (\C, 0)$,
i.e. some $\tau$ where $d \tau/ dt (0) = 0$.  See
  \cite{wall:curv} for  this fact, and for more details on the  notion of well-parameterized.

\begin{lem}
\label{wellparameter}
If $\gamma \in Germ(p)$ then $\gamma$ is well-parameterized.
\end{lem}

\paragraph{\em Proof.} Suppose not.  Then $\gamma = \sigma \circ \tau$ where $d \tau/ dt (0) = 0$
and $\sigma$ is curve germ.  We compute that  $\gamma^k = \sigma ^k \circ \tau$.
It follows that $\gamma^k$ does not immerse for any $k$, and so is not regular. Q.E.D.

\medskip

 The following proposition is key to our whole development.  In the proposition
 the statement  `$Germ(p) \sim Germ (q)$' means that for any curve $\gamma \in Germ(p)$
there is a $\sigma \in Germ(q)$ with $\gamma \sim \sigma$, and conversely.

 \begin{prop} \label{RLequivImpliesequiv} For $p, q \in \Pkn$ we have
  $Germ(p) \sim Germ (q)$
if and only if  $p \sim q$.
\end{prop}


\paragraph{\em Proof.}  Suppose that  $\gamma \in Germ(p)$ is   equivalent to $\sigma \in Germ(q)$.
Then there is a diffeomorphism $\phi \in Diff(n)$
and a reparameterization $\tau \in Diff_0 (1)$ such that $\phi \circ \gamma = \sigma \circ \tau$.
Prolonging, and using $\tau (0) = 0$, $\gamma^k (0) = p, \sigma^k (0) = q$,
 we see that $\phi^k (p) =q$.
Conversely, suppose that $p \sim q$.  Then there is a $\phi \in Diff(n)$
with $\phi^k (p) = q$.  Since $\phi^k$ preserves the distribution it preserves
the class of regular curves.  Thus, if $\gamma \in Germ(p)$ then $\phi \circ \gamma \in Germ(q)$,
showing that $\phi (Germ(p) \subset Germ(q)$.  Using $\phi^{-1}$
yields $Germ(q) \subset \phi(Germ(p))$.  Q.E.D.

\medskip

 Observe that if $\gamma^k$ is a regular curve then $\gamma^{k+1} (0)$
 is a regular point.  Consequently, if the RVT class of $\gamma$ is $\omega$
 then the RVT class of $\gamma^{k +1} (0)$ is $\omega R$.

 \begin{defn}  The R stabilization of an RVT class $w$ is
any class of the form $w R^q$, $q\geq 1$.
\end{defn}

  We restate Theorem \ref{codim1} in this R-stabilization language.
 \begin{thm}   \label{AkThm} Suppose that $p$ is in the R-stabilization of
 the class $R^k V$  and that
 $\gamma \in Germ(p)$.  Then $\gamma$ is an   $A_{2k}$ singularity.

 \end{thm}

 The proof of the theorem   requires us to develop some
    tools and is presented at the end of the next section.

 The following  is an immediate corollary of theorem \ref{AkThm} and   proposition \ref{RLequivImpliesequiv},
 \begin{cor} Each   RVT class $R^k V R^m$ (in any dimension)
consists  of a single orbit.
\end{cor}

 It is   worth pointing out a geometric  consequence of theorem \ref{AkThm}
 \begin{prop} \label{partialflag}
 A  point $p$ of type    $R^k V$  determines a unique
 partial flag  in $\C^n$ of the form  (line, plane) attached at $p_0 \in \C^n$.
 The  line is   the tangent line $p_1$ to any curve $\gamma \in Germ(p)$.
 The plane is the  tangent plane at $p_0$ to any smooth  surface germ
 containing such a $\gamma$.   \end{prop}

 The point of the proposition  is that all  $A_{2k}$ curves have embedding dimension 2,
 and  that the tangent plane at $p_0$ in the proposition  is independent of the choice of
 the particular  $A_{2k}$ curve $\gamma \in Germ(p)$.

 We are ready for the next occurring singularities.

 \subsection{Codimension 2 classes: multiplicity 3 and 4 curves.}

 The codimension two RVT classes are
precisely the $R$-stabilizations
of the classes  $$R^{s} VV , R^{s} VT; s \ge 1,$$ and $$R^{s} V R^m V; m, s \ge 1.$$
We also recall that the closure of $R^{s} VT$ intersects
$R^{s} VV$ in a class denoted by $R^{s} VL \subset R^{s} VV$ whose points
(at level $s+2$)  correspond  to those lines at level $s+1$ which lie
in the intersection of the vertical hyperplane and the critical hyperplane born from the previous level.
$R^{s} VL$ has codimension 1 within $R^{s} VV$.

Following our curve philosophy, these `next simplest' singularities in the Monster
should correspond to the ``next simplest'' curves in the classification schemes for curves.
These `next simplest' appear in Gibson-Hobbs' \cite{gib:space} work on classifying simple space curves.
Arnol'd proved that they are stably simple \cite{arn:sing}.  Arnol'd labels the corresponding
multiplicity three classes as $E_{6s+2, p, i}$ and $E_{6s, p,i}$ and various degenerations
thereof.  (See p. 23, \cite{arn:sing}.)

\begin{thm}[Classification of codimension 2 classes]\label{thm:cod2class}
 In any dimension $n$ the following holds.
 The  $R$-stabilizations of the classes $R^{s}VV\setminus R^{s} VL$
 and $R^{s} VT$   are simple,  are realized by curves of multiplicity 3,
 and their  union consists of all  points realized by
curves of multiplicity 3.
The R-stabilizations  of  $R^{s} VL$ and of
   the remaining codimension 2 classes $R^{s}VR^m V $ are
realized by curves of multiplicity 4 and are stable if and only if $s =1$.
 (In all these statements $s, m \ge 1$ are integers.)

\begin{small}
\begin{table}\label{table:table-classif-tower-simple}
 \text{Table 1: Codimension two   classes}
\renewcommand\arraystretch{1.5}
\noindent\[
\begin{array}
{|c|c|c|c|c|}
\hline
  \text{R stabilization of the RVT class} & \text Normal form   &   \text{Simple?}\\
\hline

\begin{array}{c}  R^{s} V T \\ \end{array} &  t^3 e_1 + t^{3s+1} e_2 + O(t^{3s+2})  &  Yes\\

\hline

\begin{array}{c}  R^{s}VV  \setminus R^{s+1}VL \\\end{array} &
t^3 e_1 + t^{3s+2} e_2 + O(t^{3s+4})   & Yes
\\
\hline
\begin{array}{c} R^{s} VR^m V\\ m\ge 1,  \end{array}&
t^4 e_1 + [t^{4s + 2} + t^{4s + 2m + 1}]e_2 + O(t^{4s+3})  & \text{Yes: } s =0. \text{No: } s > 0  \\ \hline

\begin{array}{c}   R^{s}VL \\\end{array} &
t^4 e_1 + t^{4s+2} e_2 + t^{4s+3} e_3 + O(t^{4s+4})  &  \text{Yes: } s =0. \text{No: } s > 0
\\
\hline

\end{array}
\]
\end{table}
\end{small}

\end{thm}


\subsection{Spatial Classes.  Bijection with Arnold-Gibson-Hobbs normal forms.}

Adding more R's to the critical classes $\omega$  of theorem
\ref{thm:cod2class} has the effect of
adding information of higher jets to the corresponding points.
The corresponding   classes  $\omega R^q$ will  break up into orbits, the stable ones
breaking up into finitely many orbits.  How many orbits? What are they?

We begin by describing those orbits which cannot be seen in the planar case.

\begin{defn}
A point p is called `spatial' if every curve in $Germ(p)$
has embedding dimension 3 or greater and at least one
of the curves has embedding dimension 3.
A point p is called `purely spatial' if every curve
in $Germ(p)$ has embedding dimension 3.
A point p is called `planar' if every curve in
$Germ(p)$ has embedding dimension 2 or more
and at least one of the curves has embedding dimension 2.
A singularity class is called spatial if all its points are spatial points.
\end{defn}

By   \ref{RLequivImpliesequiv}  a spatial point cannot be equivalent to a planar point.

\begin{prop}
 \label{1stOccurring} The 1st occuring spatial singularity classes
  occur at level 3 for $n=3$.  There are two such classes, and each is itself an orbit.
   One of these classes  forms an open subset of    $RVT$ (and so has codimension 2)
  and is  realized by the 3rd prolongations of
  those curves whose 5-jet is $(t^3, t^4, t^5)$.   The other  class is
   $RVL$  (and so has codimension 3) and is realized by the 3rd prolongations of those
   curves whose 7-jet is equivalent to   $(t^4, t^6, t^7)$.
\end{prop}

The classes of this theorem are simple stable classes, according to theorem \ref{thm:cod2class}.
It is perhaps worth noting that the class $RVT$, in any dimension greater than 2, decomposes
into precisely two orbits, one represented by the spatial curve given above, the other
represented by the planar curve $(t^3, t^4)$.

\begin{thm} \label{thm:3Dclasses} Every point of  $R^{s} VL$ is 3-dimensional.
Most of the  points of the other classes described in theorem  \ref{thm:cod2class} are spatial.
Specifically: there are dense  open subsets of $R^s VT$, of $R^s VVR$
and of $R^s V R^m V RR$ all of whose points are purely spatial.
\end{thm}

Finally,
we would like a bijection between the corresponding
R-stabilized stable simple classes and the corresponding
list of stable simple curves  of  Arnold and Gibson-Hobbs.

\begin{thm}  \label{stableClasses} For   the stable classes, $\omega = R^s VT$ and $R^s VV \setminus L$,
 there is a positive integer $q$ sufficiently large
such that the orbits of $\omega R^q$ are in bijection with the
corresponding   stable simple classes of Arnol'd's list
starting off with the appropriate Taylor series for that class,
as listed in Table 1.
\end{thm}

{\bf Example:}  $\omega = R^3 VT$.  The class begins $(t^3, t^{10}))$according to the table.
Let us use Arnol'ds notation of $(a, b+ c, d)$ to stand for curves with germ $(t^a, t^b + t^c, t^d)$.
The representative normal forms starting from $(3,10)$
are the spatial curves  $(3,10,11), (3,10 + 11, 14),  (3, 10, 14),  (3, 10+ 11, 17), (3, 10 + 14, 17),
(3,10, 17)$ and the planar curves $(3,10), (3,10 + 11),    (3, 10 + 14), (3, 10+ 17)$.
For $q \ge 19$ we are guaranteed that there are precisely 10 orbits within $\omega R^q$,
with a point in any orbit being $RL$-equivalent to one of these 10 germs.  (The value of $q = 19$
is a pessimistic  upper bound.  A value of  $q$ of about  $8$ is sufficient
to capture all 10 orbit types.)

\medskip

The reason behind this theorem  is that adding more $R$'s effectively adds   information on the derivatives
of curves.  Then, by taking $q$ large enough and fixing a point $p$ of $\omega R^q$
we have fixed enough of the Taylor series of $\gamma \in Germ(p)$
so as to be assured which one of the various stable classes it lies in.

 \section{Tools and proofs. }

We will need the following lemmas and  certain special coordinates called `KR coordinates'
after Kumpera-Ruiz \cite{kr:hist}.  These lemmas and coordinates
will also be essential tools in further sections.

\subsection{Properties of regular curves under projection}

\begin{lem}\label{prolongReg}
 Suppose that $\Gamma$ is a regular integral curve germ through $p$.
Then its  one-step projection $\pi_{k, k-1} \circ \Gamma = \Gamma_1$
 is  an  immersed integral curve.  If $p$
is not vertical then its two-step projection, $\Gamma_2 = \pi_{k, k-2} \circ \Gamma$  is   immersed.
 If $p$  is a regular point then $\Gamma_1$ is a  regular integral curve.
\end{lem}

\paragraph{\bf Proof of lemma \ref{prolongReg}.}  If $\Gamma$ is a curve germ on any manifold $Z$,
and $\pi$ is  a submersion of $Z$ onto some other manifold, then
  $d (\pi \circ \Gamma) /dt  = 0$
if and only if $d \Gamma /dt  \in ker (d \pi)$.
In our case, $\Gamma$ is immersed, and not tangent to any
critical hyperplane, so in particular, it is not tangent to the vertical hyperplane
$ker(d \pi)$ where $\pi = \pi_{k, k-1}$.  Therefore, $d \Gamma_1 /dt \ne 0$
 and $\Gamma_1$ is an  immersed curve germ.

Since $\Gamma = \Gamma_1 ^1$ we have
that $\Gamma(0) = span (d \Gamma_1 /dt (0))$,
thus, if $p = \Gamma(0)$ is not a  vertical point we have that $\Gamma_1$
is a non-vertical immersed curve, and the argument of the previous paragraph
can be repeated to yield that the two-step projection, $\Gamma_2$ is immersed.
If $p = \Gamma(0)$ is not critical (not a `C'), then $d \Gamma_1 /dt $ must
span a regular direction, so that $\Gamma_1$ is a regular curve.
Q.E.D.

\begin{lem} \label{lastVimmersed}
Suppose that  the RVT class of $p$ is $\omega_1 \omega_2 \ldots \omega_k$
and let $i \le k$ be the last occurence of the letter $V$: thus if $\omega_k = V$ then $i = k$
and if $i < k$ we have $\omega_i = V$ while $\omega_j \ne V$ for $j > i$.  Let
$\gamma \in Germ(p)$.  Then   $\gamma^{i-1}$ is
immersed and tangent to the vertical.
\end{lem}

\paragraph{\bf Proof.} Apply lemma \ref{prolongReg} iteratively until we reach level $i-1$. Q.E.D.

\begin{lem}
\label{vertmult}
Let $\gamma$ be a nonconstant analytic curve germ in $\C^n$.
Then
 $\gamma^{k+1} (0)$ is a vertical point at level $k+1$ if
and only if $mult(\gamma^{k}) < mult(\gamma^{k-1})$.
\end{lem}

The proof of this lemma requires KR coordinates and
so is postponed to the end of the next subsection.

\subsection{Kumpera-Ruiz coordinates.}

KR coordinates for the planar ($n =2$) Monster were described in detail in \cite{monz:monster}.
The generalization to general $n$   is   straightforward and detailed now. {\em For simplicity
of notation we just focus on the case $n=3$, relegating the general case to a few words near the end of
this subsection.}

We will write a KR coordinate system for $\Pp ^k (3)$ as $(x,y,z,u_1,v_1,\cdots,u_k,v_k)$. The coordinates are such that:

\begin{enumerate}
\item  $\pi_{k,j}(x,y,z,u_1,v_1,\cdots,u_k,v_k) = (x,y,z,u_1,v_1,\cdots,u_{j},v_{j})$
is the coordinate representation of the projections $\pi_{k,j}:\Pp ^k (3)\rightarrow \Pp^j (3)$, for $j\leq
k$.
\item  The last two coordinates $u_k,v_k$ are affine coordinates for the fiber.
\item There are $3^k$ KR coordinate systems covering $\Pp ^k (3)$,
corresponding to the 3 affine charts needed to cover each $\C \bP ^2$
in $\Pp^k (3) \cong \C^3 \times \C \bP ^2 \times \ldots \times \C \bP ^2$ (k times).
\end{enumerate}
We give  an inductive construction of the coordinates.
beginning with some remarks concerning  homogeneous and
affine coordinates for    projective planes.

\subsubsection{ Coordinates for a projective plane.} Suppose    the projective  plane to be
$ \bP(E)$, the projectivization of the 3 dimensional vector space $E$.
Suppose that   $E$ is endowed with   a distinguished `vertical' plane
$\Pi_{\textrm{vert}}\subset E$. Choose linear coordinates   $\theta^1,\theta^2,\theta^3$
for $E$  such that $\Pi_{\textrm{vert}} =  \{\theta^1 = 0 \}$.
The $\theta^i$ are a basis for  $E^*$ and
  $[\theta^1:\theta^2:\theta^3]$ form homogeneous coordinates on $P(E)$,
  sending a  line $\ell = span(v)  \in \bP (E)$ to the homogeneous triple
$[\theta^1 (v) : \theta^2 (v) : \theta^3 (v)]  \in \C \bP^2$.
If the line  $\ell$ is not contained in the vertical hyperplane we have $\theta^1 (v)  \ne 0$
so we may divide to get
standard  affine coordinates
$ u = \theta^2 / \theta^1, v = \theta^3 / \theta^1$
where  we use scaling to write
 $[\theta^1  :\theta^2  :\theta^3] = [1: \theta^2 / \theta^1 : \theta^3 / \theta^1 ]$.
 These coordinates cover all of $\bP(E)$ except those lines lying in $\Pi_{vert}$.

Replace  $E$ by  a rank-three distribution $D$ over a manifold $Z$.  Take the    $\theta^{i}$ to be a local    coframe for $D$. The same formulae and relations hold to yield  fiber homogeneous and fiber affine coordinates
for the prolongation $\bP D \to Z$ of $(Z, D)$.  We apply these considerations
to $(\Pp ^{k+1} (3), \Delta_{k+1})$ the prolongation of $(\Pp^k (3), \Delta_k)$.

\subsubsection{Constructing the KR-coordinates inductively}

\paragraph{{\bf The case $k=1$.}} Let  $x,y,z$ be standard coordinates on
 $\C^3$ so that  $\{dx,dy,dz\}$ form a coframe for
$\Delta_0 = T\C^3$.
Consequently $[dx: dy: dz]$ form homogeneous coordinates on $\bP (\Delta_0(x,y,z))$
and    $(x,y,z,[dx: dy: dz]) : \Pp ^1 (3) \to \C^3 \times \C \bP^2$ is
a global diffeomorphism.   There are three corresponding fiber-affine coordinates for $\Pp^1 (3)$,
depending on whether $dx \ne  0$ , $dy \ne 0$, or $dz \ne 0$.
In the case $dx \ne 0$ these   coordinates are
\begin{eqnarray}\label{eq:jet1}
u_1 = dy/dx,\\
v_1 = dz/dx   \nonumber
\end{eqnarray}
We   rewrite equations (\ref{eq:jet1}) as
\begin{eqnarray*}\label{eq:jet2}
dy - u_1 dx = 0 \\
dz - v_1 dx = 0
\end{eqnarray*}
and these two Pfaffian equations
define the distribution $\Delta_1$ on the   open set  of lines of $\Pp ^1 (3)$ for which  $dx \ne 0$.
A basis for $\Delta_1 ^*$
is formed by the restriction of $dx, du_1, dv_1$ to $\Delta_1$.

\paragraph{{\bf The case $k=2$.}}   Let  $p_2 = (p_1, \ell) \in \Pp ^2 (3)$ project onto
a point $p_1$ lying in  our  level 1    open set of   lines $\ell$
for which $dx \ne 0$.  Then
$p_1$ has KR coordinates   $(x,y,z,u_1,v_1)$.
Homogeneous coordinates
for the fibers of $\Pp^2 (3) \to \Pp^1 (3)$ are given by $[dx, du_1, dv_1]$.  The vertical hyperplane
in $\Delta_1$ is defined by $dx =0$.  We   define $KR$ coordinates $u_2,v_2$ for
a neighborhood of $p_2$
as follows:
\begin{itemize}
\item $(u_2, v_2)  = (du_1/dx , dv_1/dx)$ if $\ell$ is not vertical;
\item $(u_2, v_2) = (dx/du_1, dv_1 /du_1)$ if $\ell$ is vertical and $du_1 \ne 0$ on $\ell$;
\item $(u_2, v_2) = (dx/dv_1, du_1 /dv_1)$ if $\ell$ is vertical and  $du_1 =  0$ on $\ell$.
\end{itemize}

\paragraph{{\bf From $k$ to $k+1$} : Inductive Hypothesis.} Suppose that  KR-coordinate systems
 $\{x,y,z,u_1,v_1,\cdots,u_k,v_k\}$ have been constructed near   points $p_k \in \Pp^k (3)$,   satisfying  conditions (1) and
(2) from   the beginning of this section.
Our inductive hypothesis on the $k$-th level coordinates  is that
in each   KR coordinate system  there is  a   distinguished ordered triple of coordinates
relabeled as $(f ^k _1, f^k _2, f^k _3 )$, such that
\begin{enumerate}\label{list:previous}
\item $(df^k _1 ,df ^k _2, df ^k _3 )$ (restricted to $\Delta_{k-1}$) form a basis for $\Delta_{k-1}^*$;
\item    two of these  three coordinates $f^k _i$ are
the fiber affine coordinates $\{u_{k-1}, v_{k-1}\}$ from the previous level;
\item $df^k _1  \ne 0$ on $\ell$ where $p_k = (p_{k-1}, \ell)$  and  $\ell \subset \Delta_{k-1} (p_{k-1})$;
\item $(u_k, v_k) = ( df^k _2 /df^k _1, df^k _3 /df^k _1)$;
\item $\Delta_k$ is defined by adjoining  the Pfaffian equations
$d f^k _2 - u_k df^k _1 = 0, d f^k _3 - u_k df^k _1 = 0$
to the Pfaffian equations occurring at the lower levels $j < k$.
\end{enumerate}

Observe that under this hypothesis, a basis for $\Delta_k ^*$
is $df^k _1, du_k, dv_k$ and that the vertical hyperplane within $\Delta_k$
is defined by $df^k _1 = 0$.

\paragraph{{\bf The inductive step.}}
Take $p_{k+1} = (p_k , \ell) \in \Pp ^{k+1} (3)$
and $\{ f ^k _1, f^k _2, f^k _3 \}$ the ordered triple for $p_k$ at level $k$.
Define $\{ f^{k+1} _1, f^{k+1} _2, f^{k+1} _3 \}$ by
\begin{itemize}
\item $ (f ^k _1, f^k _2, f^k _3 )= (f^k _1, u_k , v_k)$ if $df^k _1 \ne 0 $ on $\ell$, i.e if $p_{k+1}$ is {\it not} a vertical point;
\item  $( f ^k _1, f^k _2, f^k _3 )= (u_k , f^k _1  , v_k)$ if $df^k _1 =  0 $ on $\ell$ and $du_k \ne 0$
on $\ell$;
\item $(f ^k _1, f^k _2, f^k _3 )= (v_k , f^k _1  , u_k)$ if $df^k _1 =  0 $ on $\ell$ and  $du_k =  0$
on $\ell$.
\end{itemize}
Then we have fiber-affine coordinates at level $k+1$, defined near $p_{k+1}$ by
 $$(u_{k+1}, v_{k+1} ) = (df^{k+1}_2 / df^{k+1}_1, df^{k+1}_3 / df^{k+1}_1).$$
 One checks without difficulty that $df^{k+1}_1, du_{k+1}, dv_{k+1}$
 are a basis for $\Delta_{k+1}^*$
 and that the Pfaffian system for $\Delta_{k+1}$
 is obtained by adjoining the equations
 \begin{equation*}
 d f^{k+1} _2 - u_{k+1} df^{k+1} _1 = 0, d f^{k+1} _3 - u_{k+1} df^{k+1} _1 = 0
 \end{equation*}
 to the equations arising from the lower levels.
 The inductive hypothesis for the $k$-th step implies the hypothesis for the
 $(k+1)$-th step.  We have completed the construction of the KR-coordinate systems.\\

\begin{eg}\label{jetEg}{\bf At Cartan points KR coordinates are jet coordinates.}

Take a Cartan point $p \in C^k (3) \subset \Pp^k (3)$. Let   $p_1 \in \Pp^1 (3)$ be its
projection to level $1$, a line in $\C^3$. Choose linear coordinates $x, y, z$
on $\C^3$ such that    $dx \ne 0$ on   $p_1$.
At each step $j \le k$ of the KR construction we have $f^j _1 = x$,
so that  $u_k = du_{k-1}/dx, v_k = dv_{k-1}/dx$.
Combining equations we get
$$(u_j , v_j) = ( d^j y /dx ^j,  d^j z /dx^j)$$
These are standard jet coordinates for maps from $\C$ to $\C^2$. The Pfaffian system
for $\Delta_k$ is $du_{j-1} - u_j dx = 0, dv_{j-1} - v_j dx = 0$,
$j = 1, \ldots , k$ which is the standard distribution on
the jet space $J^k (\C, \C^2)$.

In the case of  $p \in C^k (n)$ for general dimension $n$
the computations are nearly   identical.  Take
 linear coordinates $x_1, \ldots , x_n$  on $\C^n$ for which $dx_1 \ne 0$
 on the line of $p_1$. Set  $x = x_1$ and work over the open set at  level 1 of those lines for which
$dx \ne 0$.  At level $j$  the KR-fiber coordinates
$u_{2, j}  , \ldots, u_{n-1, j} $ satisfy  $u_{i, j}  = d^j (x_i) /dx ^j$ in a \nbhd of any curve in $Germ(p)$.
These are standard jet coordinates for maps from $\C$ to $\C^{n-1}$. The Pfaffian system
for $\Delta_k$ near $p$ is  given by $du_{i, j-1} - u_{i,j} dx = 0$,
$j = 1, \ldots , k$, and $i = 2, \ldots , n-1$.  These define the   standard distribution on
the jet space $J^k (\C, \C^{n-1})$.

\end{eg}

\medskip

\paragraph{\bf Finishing the Proof of theorem \ref{Cartan} .}

Proposition \ref{CartanNV} establishes all assertions of the theorem except that concerning
the identification of a \nbhd of a Cartan point with a \nbhd of $J^k (\C^1, \C^{n-1})$.
This assertion follows directly from the  example \ref{jetEg} above.  Q.E.D.

\medskip

The construction for KR coordinates on the Monster
$\Pp^k (n)$ for general $n$  proceeds in a nearly identical manner.
The  main    difficulty is {\em notational}.
Label  KR coordinates for $\Pp ^k (n)$ as
$(x_1, \ldots, x_n, u_1 ^1, \ldots , u_{n-1} ^1, \ldots ,u_1 ^k, \ldots , u_{n-1} ^k)$.
The affine fiber coordinates $u_1 ^k, \ldots , u_{n-1} ^k$ are
built out of the previous   level according to  $u_j ^k = d f_j  / dg$
where $f_1, \ldots , f_{n-1}, g$ an n-tuple of coordinates    selected from among
$\{ u_1 ^{k-1}, \ldots , u_{n-1}^{k-1}  \}$ and one of the coordinates coming from a level less than $k-1$.
These coordinates are such that  the  $df_i$ together with $dg$ form a basis for
$\Delta_{k-1}^*$, and  $dg \ne 0$ on  the line  $\ell$ of
of $p = (p_{k-1}, \ell)$   which the coordinate \nbhd must cover.
We leave further  details  to the reader.

 \medskip

\paragraph{\bf Proof of Lemma \ref{vertmult}.}  The point $\gamma^{k+1} (0)$
is vertical if and only if the curve $\gamma^k$ is tangent to the vertical
at level $k$. In KR coordinates, $\gamma^{k}$ is represented  by adjoining
$n-1$ new fiber
affine coordinates   to the KR  coordinate representation of $\gamma^{k-1}$.  The curve
    $\gamma^k$ is tangent to the fiber, i.e. to the vertical,
if and only if the order  of at least one of the new  fiber coordinates
is less than the orders of {\it all}    of   the previous  coordinates,
those coordinates representing $\gamma^{k-1}$.  In other words,
if and only if $mult(\gamma^k) < mult(\gamma^{k-1}).$ Q.E.D.

\medskip

\begin{eg}\label{A2keg}{\bf   The $A_{2k}$ singularity.}
Start with the curve $c$   given by  $x = t^2, y =t ^{2k +1}$.
Then $u_1 = dy /dx = \frac{2k+1} {2} t^{2k-1}$ defines the first prolongation of the curve in KR coordinates
$(x,y, u_1)$.
The $j$th prolongation, $j \le k$  is given by adding fiber coordinate
  $u_j = d^j y /dx^j = c_j t^{2(k-j) +1}$ to the previous $u_i$,  $i < j$, where  the
$c_j$ are positive rational numbers.  Consequently, referring to lemma  \ref{vertmult}
and the spelling rules, we see that the first $k$ letters of the RVT code for $\gamma^N (0)$,
$N \ge k$, are R's.  At level $k$, the curve becomes immersed, tangent to the vertical,
with lowest order coordinate $u_k$ having order $1 < 2$.
It follows that the $(k+1)$-st letter of the RVT code is $V$.
 At level $k+1$ we compute that
the new KR coordinate is $dx/du_k = c_{k+1} t $ representing a regular direction, since $c_{k+1} \ne 0$.
The curve $c(t)$ regularizes at level $k+1$.  Now proposition \ref{1regularizing} yields that the code of $c^{k+1 + s} (0)$
is   $R^k V R^s$.

\end{eg}

\subsection{Preparing Curves.}   Given a curve germ $\gamma(t)$ in $\C^n$ we can always, by linear change of coordinates,
find coordinates $x_i$  centered at $\gamma(0)$ so that when the curve is expressed  in these coordinates
we have
$$ord(x_i (t)) <  ord(x_{i+1} (t)).$$
  In these coordinates
$ord(x_1 (t)) = m$ where $m = mult(\gamma)$.
Finally, we can reparameterize the curve
so that $x_1 (t) = t^m$.
When such coordinates and a parameterization  are chosen, we will say we have {\it prepared } $\gamma$.
Thus, a prepared curve, is given in these coordinates by
$\gamma(t) = (t^m, x_2 (t), x_3 (t), \ldots )$ with the $x_a (t)$
power series in $t$:
$x_a (t) = \Sigma_{j> m_a} A_{a , j}  t^j$, $a =2, \ldots, n$ and
$m < m_1 < m_2 < \ldots $ etc.

\subsection{Proof of  the  theorem  \ref{AkThm}, the ``$A_{2k}$ theorem''.}

By    example \ref{A2keg} immediately above,   the $A_{2k}$ singularity realizes a point of   type $R^{k} V R^s$
upon   $k+1 +s$  prolongations.
To finish the proof, we must show that  if $q \in R^{k} V R^s$,
and if $\gamma \in Germ(q)$ then $\gamma$ is  an   $A_{2k}$ singularity.

 It follows from lemma \ref{prolongReg} that if  $q$ lies one step over $p$ and is a regular point,
then $Germ(q) \subset Germ(p)$.  Iterating, we see that it
 suffices to show that if $p \in R^{k} V$
and $\gamma \in Germ(p)$ then $\gamma$ has type $A_{2k}$.
It suffices in turn to show
that for all  $p \in R^{k} V$ and all $\gamma \in Germ(p)$,
we have   $mult(\gamma) = 2$.
{\em This is sufficient because  lemma \ref{wellparameter} tells us
that $\gamma$ is well-parameterized, and  any well-parameterized curve $\gamma$ of multiplicity $2$
is an   $A_{2j}$ singularity for some $j$.}   Finally,  we must
have $j = k$ since the RVT code of  our point $p = \gamma^{k+1}(0)$ is $R^k V$
while the RVT code of a point $\gamma^{k+1} (0)$ for $\gamma$ an $A_{2j}$ singularity  is $R^{k+1}$ if $j > k$
and     $R^j V R^{k-j}$ if $j < k$,  according to  the computation at the end of the last section.

Write $m = mult(\gamma)$. By  definition,   $\gamma^{k+1}$ is a regular integral curve.
By lemma \ref{prolongReg}, its one-step projection $\gamma^{k}$ is an immersed
curve. This immersed curve is tangent to the vertical space,
since $\gamma^{k+1} (0)$ is a vertical point.     Thus,
in a KR coordinate system one of the fiber coordinates, say $u$,
for $\gamma^k$ has Taylor series   $at + \ldots$, $a \ne 0$.
Since the $\gamma^j (0)$, $j < k+1$ are all regular points,
we have, by lemma \ref{vertmult},  that $mult(\gamma^j ) = m$ for all these $j < k+1$.
 As discussed just above,  we can ``prepare'' our curve, which is to say
choose coordinates $x = x_1, y = x_2, x_3 \ldots, x_n$
so that when $\gamma(t)$ is expressed in these coordinates
we have $ord(x_i (t)) < ord(x_{i+1}(t))$
and $m =ord(x(t))$.   Write $m_i = ord(x_i (t))$.
The corresponding KR  fiber-affine coordinates along  the first prolongation of $\gamma$
are  $u_i = dx_i /dx$ while the subsequent fiber-affine KR coordinates at level $j$, $j < k+1$
are of the form $d^j x_i /dx^j$ and so have    order $m_i - j\mbox{ }m$.
 (Note $\gamma^j (0)$ is a Cartan point for $j \le  k$.  See the example above on Cartan points.)
It follows that    $mult(x_i) > k \mbox{ }m$ for $i \ne 1$.
 Since     $\gamma^{k}$ is immersed
  and tangent to the vertical by lemma \ref{lastVimmersed}   the coordinate
  with the next smallest order after $x(t)$, namely $y =x_2$,   must have   multiplicity $k\mbox{ }m +1$,
to yield   order $1 = (k\mbox{ }m+1) -k\mbox{ }m$ for one of the fiber coordinates at  level $k$.
 Write this     affine coordinate   as $u = d^k y /dx^k$.
The   KR fiber coordinates  at   level $k+1$ are then
$U = dx/du$ and $V_j = dw_j/ du$ where $w_j$ represent the rest of the   affine coordinates at level $k$
($w_j =  d^k x_j /dx^k$, $j > 2$).

To finish the proof, we will need to show that
the  regularity of $\gamma^{k+1}$ implies that $ord(U(t)) = 1$.
To establish this fact  we   need to express  the  critical planes through $p$ in KR coordinates.
A basis for $\Delta_{k+1}^* (p)$ is $du, dU$ and the $dV_j$.   $U$ and the $V_j$ coordinatize the fiber through $p$
so that the vertical hyperplane through $p$ is defined by   $du = 0$. There is one more critical
hyperplane, $\delta_k ^1 (p)$,  through $p$ and this hyperplane arises   from  the baby monster through $p_k$,  one level down.  At level $k$
the fiber coordinates are $u$ and the $w_j$, hence any curve lying in this fiber has
$x = const.$, and so $dx =0$.  Since $U = dx/du$ we must have that $dU = 0$ along   prolongations of   curves lying in
the fiber through $p_k$, showing that $\delta_k ^1 (p)$ is given by $dU = 0$.
Thus,  the regularity of   $\gamma^{k+1}$ is equivalent to  $du \ne 0$ and $dU \ne 0$ along $\gamma^{k+1}$.
The second condition implies that $mult(U) =1$.  But $ord(U) = ord(x) - ord(u) = m-1$,
implying that
$m = 2$, which is to say,  $mult(\gamma) =2$.  Q.E.D.

\medskip
{\bf Proof of proposition \ref{partialflag}. }
Fix a point $p \in R^k V$ and a curve $\gamma \in Germ(p)$. As per the
above proof, there are   coordinates  $x, y, z, \ldots, $ on $\C^n$ associated to $\gamma$ such that
$\gamma = (t^2, t^{2k+1}, 0, \ldots ) + O(t^{2k+2})$.  Moreover, by expressing prolongations
of curves near $\gamma$ in  the associated KR
coordinates, as per the proof of proposition,  we see that {\it any} $\bar \gamma \in Germ(p)$, is,
after  reparameterziation,
given by  $\bar \gamma = (t^2, a t^{2k+1}, 0, \ldots) + O(t^{2k+2})$, with $a \ne 0$.
Any surface containing any one of these curves $\bar \gamma$ agrees with the surface defined by the equations
$z = 0, \ldots , x_n = 0$ up to terms in $x, y$ of the form $xy,  x^{2k}$. The tangent space to any such
surface is the $x,y$ -- plane.
Q.E.D.

\subsection{Proof of the codimension 2 theorem, theorem   \ref{thm:cod2class}.}

There are only finitely many RVT classes at any given level,
and   they exhaust the monster tower at that level.    Thus
to show a given RVT class $\omega$  is simple it suffices to show that
\begin{itemize}
\item[(a)] $\omega$ consists  of finitely many orbits and
\item[(b)] every class adjacent to $\omega$,  i.e. whose closure contains $\omega$, consists of finitely many orbits.
\end{itemize}
The RVT classes of codimension 2 and length $s+ m + 1$
are adjacent to the classes $R^s V R^m$
and $R^{s + m + 1}$.  (See proposition \ref{RC}.)  We have seen (Theorem \ref{Cartan}, Proposition \ref{CartanNV},
and Theorem \ref{AkThm}) that these classes
each consist of a single orbit,  so criterion (b) holds.  It remains to check criterion (a).
{\bf We argue class by class}, following the structure of the   proof of the $A_{2k}$ theorem \ref{AkThm}.
The argument for every class $\omega$ of the theorem begins in
the same way.

Suppose that $\omega$ is one of the codimension two RVT classes,
that $p \in \omega$ and that $\gamma \in Germ(p)$.
Since $\omega$ is not a codimension 0 or 1 class, we know that
$mult(\gamma):= m  > 2$.  We can choose coordinates $(x,y,z, \ldots)$ on
$\C^n$ and a parameterization
 so that our curve has the prepared form $\gamma = (x(t), y(t), z(t), \ldots)$ :
  $x = t^m, y = t^{m_2} + o(t^{m_2}),  z = \epsilon  t^{m_3} + O(t^{m_3 + 1}), \ldots$
with  $m < m_2 < m_3 < \ldots$  and $\epsilon =1$ or $0$.   Since the first $s$ letters of $\omega$
are $R$'s we know from lemma \ref{vertmult} and the spelling rules
  that  the multiplicity of the curves
$\gamma, \gamma^1, \ldots, \gamma^{s-1}$
are   all  equal to $m$, and that  their $KR$-coordinates up to level $s$
are of the  Cartan form above, as expressed in   example  \ref{jetEg}.
In particular, at level $s$ we have fiber affine coordinates
$$u =  d^{s} y /dx^{s} , v = d^{s} z /dx^{s}, \ldots $$
Since $\omega_{s+1} = V$ we have, by lemma \ref{vertmult},  that
$ord(u) < m = ord(x)$, and because $m_2 < m_3 < \ldots$ etc.
that $ord(u) < ord(v) < \ldots$ etc.  At order $s+1$  the new fiber affine coordinates
are
$$u_2 = dx/du,  v_2 = dv/du, \ldots \mbox{ etc.}$$
A basis for $\Delta_{s+1} ^*(p_{s+1})$ is $du, du_2, dv_2, \ldots$ etc.

There are exactly two  critical planes at level $s+1$ and these are given by
 \begin{itemize}
 \item $du = 0$ (the vertical hyperplane) and
 \item $du_2 = 0$ (the tangency hyperplane $\delta_s ^1$).
 \end{itemize}
\medskip

\paragraph{\bf The class $R^{s} VT$.}   Let   $p$ be a point in this class and  $\gamma \in Germ(p)$.
 No letter past the $(s+1)$-th in the   code is a
$V$, so lemma
\ref{lastVimmersed}
implies that   $\gamma^{s}$ is immersed
and tangent to the fiber.    Thus, in the notation above,   $u =  d^{s} y /dx^{s} $ has order $1$ while
the other fiber coordinates at this level have order   greater than or equal to  $2$.
 Since our  point represents $R^{s} VT$ the immersed curve $\gamma^{s+1}$
  must be tangent to the T hyperplane, so that $du_2 = 0$ along the tangent to $\gamma^{s+1}$.   This equality  simply
  asserts that $ord(x) > 2$ since $ord(u_2) = ord(x) -ord(u) = ord(x) - 1 > 1$. At level $s+2$
  the new fiber coordinates are
  $$u_3 = du_2/ du = d^2 x /du^2,  v_3 = dv_2 / du = d^2 v /du^2, \ldots.$$
  A basis for $\Delta_{s + 2}^*$ is $du, du_3, dv_3, \ldots$.
  The tangency hyperplane $\delta_s ^2$  through our level $s+2$ point
  which corresponds  to the baby monster issuing from  level $s$
  is defined by $du_3 = 0$ while the vertical hyperplane through that point is defined by $du=0$ again.
  Since $\gamma^{s+2}$ is a regular curve, we must have $du_3 \ne 0$ along $\gamma$ at $t =0$.
  Since $u_3 = du_2/ du = d^2 x /du^2$ has order $m-2$ we get
  that $m-2 =1$ so that $m =3$.

  We now have that $\gamma = (t^3, t^{m_2},  \ldots)$.
  From the expression  $u = d^{s} y /dx^{s}$ we have that
  $ord(u) = m_2 - (3s)$.  But $ord(u) =1$ so we must have
  $m_2 = 3s+1$.    Thus $\gamma = (t^3, t^{3s+1}, \ldots) + o(t^{3s+1})$.
Results from the  singularity theory of
  curves (see e.g. \cite{arn:sing} )   asserts that    any   curve whose $(3s+1)$-jet has  this
  form   is simple and
  is RL equivalent to one of
  a finite number of curves   listed as $E_{6s, p, i}$ in  Arnold (op.cit.).
  All   curves in this list are 2 or  3-dimensional.

    We have    established that $\gamma$ is simple. By Proposition \ref{RLequivImpliesequiv} , and  the adjacency argument (a) above,    $p$ is simple.   This finishes the proof for
    the case $R^s VT$ and its R-stabilizations.
     \medskip

\paragraph{\bf The class {$R^{s}VVR^m\setminus R^{s}VL R^m$.}}
By lemma \ref{prolongReg}  if $p$ is a point in this class and  $\gamma \in Germ(p)$ then $\gamma$  must immerse at level $s+1$ and is tangent to
the vertical at that level.
As discussed above, a  basis for $\Delta^*_{s+1}$  is  $\{du,du_2,dv_2, \ldots \}$
and in these linear coordinates,
the vertical hyperplane within $\Delta_{s+1}$  is given by $du = 0$ while the tangency (non-vertical  critical )hyperplane    is given by $du_2=0$.    Since   $\gamma^{s+2}(0)$ is a V point, we have
that, at $t =0$, $du = 0$  and since
 $\gamma^{s+2}(0)$  is  {\it not} a T point (i.e. not an   L) point, and $\gamma^{s+1}$ is
 immersed, we have that  $du_2 \ne 0$.   It follows that  $ord(u)>1$
while $ord(u_2) =1$.

At level $s+2$ we have new KR coordinates
$$u_3 = du/du_2,  v_3 = dv_2/ du_2, \ldots$$
and  $\{ du_2,du_3,dv_3, \ldots  \}$ form a basis for $\Delta^{*}_{k+2}.$
 In this basis the vertical hyperplane is given by $du_2 = 0$
 while the  (non-vertical) critical hyperplane $\delta_{s+1}^1$
 is given by  $du_3 = 0 $. Since $\gamma^{s+2}$ is regular, we   have $du\neq 0,du_3\neq 0$ and so  $u_2,u_3$   both have order $1$.
But $ord(u_2) = ord(x) - ord(u)$ and $ord(u_3) = ord(u) - ord(u_2) = 2 ord(u) - ord(x)$. Thus
$$
\begin{array}{c}
 ord(x) - ord(u) = 1   \\
 2 ord(u) - ord(x) = 1
\end{array}
$$
The unique solution to this linear system is  $ord(x)=3, ord(u)=2$.
 Since $ord(u) = ord(y) - (s)ord(x)$ we get   $y = t^{3s+ 2} + O(t^{3s+3}).$
 We now have $\gamma$ in the desired normal form,
 $(t^3,  t^{3s + 2}, \ldots , 0) + O(t^{3(s+1)})$.

 Results from the  singularity theory of
  curves (op. cit. see eg  \cite{arn:sing} )   asserts that    any   curve whose $(3s+2)$-jet has  this
  form   is simple and
  is RL equivalent to one of
  the finite number of curves   labeled $E_{6s+2, p, i}$ by Arnol'd.  See also Gibson-Hobbs \cite{gib:space}.

Again,  we have    established that $\gamma$ is simple, and by the adjacency arguments that $p$ is simple.
    We have established that $\gamma$ has multiplicity $3$ and is either a
    planar or a  strictly 3-dimensional curve.    This finishes the proof for
    the case $R^s VV \setminus R^s VL$ and its R-stabilizations. \\

\medskip

\paragraph{\bf The class $R^{s} V L$.}
We continue with the notation and coordinates from the previous case.
At level $s+1$, $\gamma^{s+2}$ is immersed by lemma \ref{prolongReg} and
 $\Delta^{*}_{k+2}$
has basis $\{du,du_2,dv_2, \ldots  \}$, but now  $du=du_2=0$
along the curve $\gamma^{s+1}$ since the letter `L'  of the code asserts that the curve is tangent to a line
lying in both  critical planes.
It follows that both $u$ and $u_2$
have multiplicity $2$ or greater, while  the coordinate of lowest order,
$v_2$, must have order 1 since the curve is immersed.

At    level $s+2$,   we   have fiber affine coordinates:
\begin{equation*}
u_3 = \frac{du}{dv_2} , v_3 = \frac{du_2}{dv_2},
\end{equation*}
and   basis $\{ dv_2,du_3,dv_3, \ldots  \}$ for $\Delta^{*}_{s+2}$. Through `L' points
there are (at least) 3 critical hyperplanes, namely  $dv_2 = 0$,   $du_3 = 0$ and   $dv_3=0$.  (When $n =3$
there are exactly 3 such planes.) More generally, for any $n$ if $\omega$ is of the
form   $\alpha RVL$ where $\alpha$ is arbitrary and $p \in \omega$ then  there are exactly three critical planes
through $p$.   See the section below on intersection combinatorics of critical planes.)
Since $\gamma^{s+2}$ is regular  all three of $dv_2, du_3$ and $dv_3$
must be nonzero along the tangent to   $\gamma^{s+2}$.
That is, $ord(v_3) = ord(u_3) = 1$.
Using  $ord(v_2) =1$  and the relations defining $v_2$,
$u_3, v_3$     we derive,
\begin{itemize}
\item $ord(v_2) = 1 = ord(v) - ord(u)$,
\item $ord(u_3) = 1 = ord(u) - ord(v_2) \Rightarrow ord(u) = 2,$
\item $ord(v_3) = 1 = ord(u_2) - ord(v_2) \Rightarrow ord(u_2) = 2.$
\end{itemize}
Putting  the results of the second and third equations
into the first, we deduce that $ord(v) = 3$. Since $ord(u_2) = ord(x) - ord(u)$ we obtain $ord(x) = 4$ and from $ord(u) = ord(y) - s \mbox{ }ord(x)$ we derive $ord(y) = 4s + 2$.
From $v = d^{s} z / dx^{s}$ we see that $ord(v) = ord(z) - s \mbox{ }ord(x)$
so that $ord(z) = 4s + 3$.  Scaling now,
we can put our curve into the form
$(t^4, t^{4s+ 2}, t^{4s+ 3}, \ldots ) + (0, O(t^{4s + 3}) , O(t^{4s+ 4}, \ldots)$.
A diffeomorphism of the form $(x, y, z, \ldots) \to (x, y + a xz, z)$ kills the middle
 $O(t^{4s  + 3})$ and puts our curve
 into the desired form.

 According to Gibson-Hobbs, or Arnol'd,
 the curve above is simple if and only if $s =1$.    This finishes the proof for
    the case  $R^s VL$.

\medskip

\paragraph{\bf   Last case:  $R^s V R^m V$ .}

Fiber affine coordinates at level $s$ are $u, v, \ldots$ with
\begin{equation*}
u = \frac{d^{s}y}{d x^{s}} , v = \frac{d^{s}z}{d x^{s}} , \ldots
\end{equation*}
A basis for  $\Delta^*_{s} $ is $dx, du, dv, , \ldots$ and in these
linear coordinates on  $\Delta_{s} $ the vertical hyperplane is  $dx=0$.
Since the $(s+1)$-th letter in the code is a $V$, $\gamma^s$ is tangent to the vertical:
 $dx(\gamma^{s})'(0)=0$.
 Since $ord(y) < ord(z)$ we have that the coordinate with lowest order at   level $s$
 is $u$.
At level $s+1$     fiber coordinates are:
\begin{equation*}
u_2 = \frac{dx}{du}, v_2 = \frac{dv}{du},
\end{equation*}
  $\Delta^{*}_{s+1}$ has basis  $du,du_2,dv_2 $
  and in these
linear coordinates for  $\Delta_{s+1} $ the vertical hyperplane is  $du=0$ while the
critical hyperplane arising from the baby-Monster one level down is  given by $du_2 =0$.
Since $\gamma^{(s+1)}$ is tangent to a regular direction,
we have
\begin{equation}
 \label{orders}
ord(u) = ord(u_2) < ord(v_2)
\end{equation}
 From level $(s+2)$,   all the way up to level $s+m +1$ the dominant (lowest order) coordinate
continues to be  $u$  and the subsequent fiber coordinates are derivatives with respect to $u$:
\begin{eqnarray*}
u_3 = \frac{du_2}{du}, v_3 = \frac{dv_2}{du},
\ldots, u_{m+2} = \frac{du_{m+1}}{du}, v_{m+2} = \frac{dv_{m+1}}{du}  \\
\end{eqnarray*}
and $ord(u_i) < ord(v_i) <  \ldots$ for $i = 3, 4, \ldots , m+2$.
We have that $\gamma^{s+m+1}$ is vertical
due to the occurrence of the final V in the RVT code and
so   $u$ is no longer the dominant coordinate at this level: $ord(u) > min\{ ord(u_{m+2}), ord(v_{m+2})\}$.
 Lemma \ref{prolongReg}  implies that   $\gamma^{s+m+1}$ is immersed so that   $u_{m+2}$ is now the dominant coordinate, with order $1$.

At level $(s+m+2)$
\begin{eqnarray*}
u_{m+3} = \frac{du}{du_{m+2}}, v_{m+3} = \frac{dv_{m+2}}{du_{m+2}},\\
\Delta^*_{s+m+2} = span \langle du_{m+2},du_{m+3},dv_{m+3, }\rangle.
\end{eqnarray*}
 the {\em vertical hyperplane} is  $du_{m+2} = 0$ and the critical hyperplane arising from the   baby Monster one level down
 is   $du_{m+3}=0 $. Since $\gamma^{s+m+2}$ is regular  we   have $du_{m+2}, du_{m+3}\neq 0$ along $\gamma^{s+m+2}$. Therefore,  the order of  $u_{m+3}$ is also 1.   But
\begin{equation*}
1 = ord(u_{m+3}) = ord(u) - ord(u_{m+2}) = ord(u) - 1 \Rightarrow ord(u) =2.
\end{equation*}
Moreover, from equation \ref{orders} and $u_2 = dx/du$  we deduce $2 =ord(u_2) = ord(x) - ord(u) \Rightarrow ord(x) = 4$. \\
In conclusion:
  \begin{eqnarray*}
  ord(u) = ord(y) - s\mbox{ } ord(x)\Rightarrow ord(y) = 4 s + 2.
   \end{eqnarray*}
 Finally, $ord(z) < ord(y)$.

Knowing that $x = t^4 + \cdots$ and $y = t^{4s+2} + \cdots$ we shall determine what is the smallest non-vanishing term in $(t^4,t^{4s+2},0) + o(4s+2)$ which makes this germ well-parameterized.  \\
Let us say the first non-vanishing in $y$ is $t^N$ for $N > 4s + 2$, which allows us to write $$y(t) = t^{4s+2} + d_0 t^N + \cdots.$$

By the same arguments as in the beginning of the proof, we have

$$u(t) = \frac{d^s y}{dx^s} =  c_1 t^2 + d_1 t^{N - 4s} + \cdots$$
and $d_1$ is proportional to the original constant $d$ in the definition of $y$.

Now either
\begin{itemize}
\item $ord(v) < ord(x)$ or
\item $ord(v) \geq ord(x)$.
\end{itemize}

In the former case, we have $ord(v_2) < ord(u_2)$ contradicting equation \ref{orders} above. Therefore $ord(v)\geq ord(x)$. Differentiating this inequality yields $ord(v_i)>ord(u_i)$ up to level $s+m+2$ which justifies why we can safely ignore the `$z(t)$' and other component terms `$(x_i(t),i>3)$' in the demonstration.

We compute

$$u_2(t) = \frac{dx}{du} =\frac{t^4}{c_1 t^2 + d_1 t^{N - 4s} + \cdots} = c_2 t^2 + d_2 t^{N - 4s} + \cdots.$$ It is not hard to convince ourselves that the coefficients $c_2,d_2$ in this last equation are rational functions of previous Taylor coefficients.
Likewise,
$$u_3(t) =\frac{du_2}{du} = c_4 + d_4 t^{N - 4s - 2} + \cdots$$ and moreover $c_4$ and $d_4$ are both non-zero. Continuing in this fashion, we can demonstrate that
$$u_{m+2}(t) = c_{m+2}t^{N - 4s - 2m} + \cdots,$$ and $m\geq 1$.
From our earlier logic, $u_{m+3} = \frac{du}{du_{m+2}}$ has a non-zero linear part as a function of $t$, and therefore from
$$u_{m+3}(t)  = \frac{d(c_1 t^2 + d_1 t^{N - 4s - 2} + \cdots)}{d( c_{m+2}t^{N - 4s - 2m} + \cdots)},$$
we can conclude, after deriving the resulting quotient series that
$$2 - (N - 4s - 2m)=1\Rightarrow N = 4s + 2m + 1.$$
Notice that this latter power is {\em odd},    guaranteeing that our initial germ is well-parameterized. The parameter $d_0$ can be rescaled to 1 by RL equivalence and therefore does not correspond to moduli. Also, had we assumed the non-vanishing term $t^N$ appeared in the $z$-component the resulting curve would be spatial. \\

The corresponding germs are simple if and only if  $s=1$.
We refer to the lists  in Gibson-Hobbs and Arnol'd.
    This finishes the proof for
    the last case  $R^s V R^m V$ and its R stabilizations.  Q.E.D.

\subsection{Proofs regarding Spatial classes.}

  {\bf Proof of \ref{1stOccurring}.}

  The case of $RVL$ has already been proved.

  We proceed to  the $RVT$ case.   Any point $p_3$ of this type
  projects one level down to a point of type RV.
  Fix a point $p_2$ at level $2$ of type RV.  All such points are equivalent.
  As per proposition \ref{partialflag}, $p_2$ determines a
  partial flag passing through $p_0$.    We choose coordinates $(x,y,z)$  so that $p_0$
  is the origin, the line ($p_1$) is the $x$-axis
  and the plane is $x,y$ -- plane.  Now, consider  the locus of all
  points $p_3$ of type RVT lying over $p_2$, and consider
  the set of all resulting curves $\gamma \in Germ(p_3)$.
    According to the KR computations done above
    ( the case $R^s VT$ with $s = 1$ within the Proof of the codimension 2 theorem, theorem   \ref{thm:cod2class}),
    and with the base coordinates
  being these $(x,y,z)$, any such curve $\gamma$
  has the form
  $$\gamma = (t^3, a_4 t^4 + a_5 t^5,  b_5 t^5) + O(t^6),  a_4 \ne 0 $$
  after a reparameterization.

We compute the corresponding fiber  KR coordinates

$$\text{level 1}:  u =d y/dx = {4 \over 3} a_4 t + {5 \over 3} a_5 t^2 + O(t^3),v = dz/dx = {4 \over 3} b_5 t  + O(t^3);$$

$$\text{level 2}:  u_2 =d x/du = {9  \over {4 a_4}} t ^2 + O(t^3),v_2 = dv/du = {5 \over {2 a_4}} b_5 t  + O(t^2);$$

$$\text{level 3}:  u_3 =d u_2 /du = c t  + O(t^2 ),v_3 = dv_2/du = {15 \over {8 a_4}} b_5  + O(t).$$

In these coordinates, the 3rd prolongation of $\gamma$ is $(0,0,0; 0,0; 0, 0; 0, {15 \over {8 a_4}} b_5 )$.
Now $\gamma$
is spatial if  $b_5 \ne 0$. Said invariantly: every curve sharing 5-jet
with $\gamma$ is spatial if and only if  $b_5 \ne 0$.  We have proved that those points with
$b_5 \ne 0$ are purely spatial.   Now any  curve $\gamma$  with 5-jet of the given form, having  $b_5 \ne 0$
is RL equivalent to $(t^3, t^4, t^5)$.  The proposition is proved.
Q.E.D.

\smallskip

Note that the curves with $b_5 = 0$ include the planar curve $(t^3, t^4, 0)$.   Their 3rd prolongations forms
the rest of the class $RVT$ -- the set of points  having planar curves in their germ.

 \medskip

\subsection{Proof of Theorem \ref{thm:3Dclasses}.}
\label{RsVT}

Looking back at the   previous proof, that of proposition  \ref{1stOccurring},  we observe that the
determining factor was whether or not the coefficient  $b_5$ in $\gamma$'s
Taylor expansion vanished.  The point   at level 3 in RVT
was spatial if and only if the  coefficient $b_5 \ne 0$.
This coefficient occurs in the 2-jet of the curves $\gamma^1$ at level
1.  It determines the outcome of points at level $3 = 1 +2$.

To clarify our understanding
we use the following lemma.
In the lemma we use the symbol $=_{rep}$ to mean ``equal up to reparameterization''.

\begin{lem}\label{lem:eq} Let  $\Gamma,\tilde{\Gamma}:(\R,0)\rightarrow (\Pkn,p)$  be immersed integral curve germs
passing through $p$ and    $q\geq 1$ an integer.  Then    $\Gamma^{q}(0) = \tilde{\Gamma}^{q}(0)$
if and only if $j^{q} \Gamma =_{rep}  j^{q} \tilde{\Gamma}$.

The corresponding   map $j^q \Gamma \mapsto \Gamma^q (0)$
from q-jets of immersed  integral curves through
$p$ to points in the Monster lying q steps over $p$ is algebraic.
\end{lem}

 The lemma is a generalization from
 $n =2$ to general $n$ of a   lemma crucial to  the book
(\cite{monz:monster} , pp. 54).
We give an   conceptual proof
as an alternative to  the book's coordinate-based proof.

\paragraph{\bf Proof of Lemma \ref{lem:eq}}.
 Forget the distribution $\Delta_k$ on  $\Pkn$ for a moment, treating
$\Pkn$ as a  complex manifold $Z$.   The lemma is certainly true
in this more relaxed situation: for immersed curves germs $\Gamma,\tilde{\Gamma}$ in $Z$, and for $q\geq 1$
we have $\Gamma^{q}(0) = \tilde{\Gamma}^{q}(0)$
if and only if $j^{q} \Gamma =_{rep}  j^{q} \tilde{\Gamma}$.  Here,
the curves  $\Gamma^q,  \tilde {\Gamma}^q$ are curves in the q-fold prolongation $\Pp^q (Z)$ of $(Z, TZ)$.
The lemma  is true in this relaxed situation because  space $\Pp^q (Z)$ with its distribution
is  locally isomorphic to the prolongation tower $\Pp^q (N)$, $N = dim(Z)$ with its $\Delta_q$ and the points $\Gamma^q (0), \tilde \Gamma^q (0)$ are Cartan points in $\Pp^q (Z)$.Now  use the KR    computations as per
 example \ref{jetEg} which relate   neighborhoods  of
Cartan points to neighborhoods in the appropriate jet  space of curves, thus establishing the relaxed version of the lemma,
including the algebraic nature of the map.
To finish the proof of the original lemma,  simply observe that
$\Pp^{k+q} (n) \subset \Pp^q (Z)$ is an algebraic  submanifold,  the curves $\Gamma^q, \tilde \Gamma^q$ lie in this
submanifold, and their q-jets form an algebraic submanifold of all q-jets of immersed curves through $p$. Q.E.D.

\medskip

\paragraph{\bf Case $R^s VT$.}  \label{RsVT} We can almost copy the  previous proof.
Fix a point $p_{s+1}$ representing $R^s V$.
 All such points are equivalent.
  As per proposition \ref{partialflag}, $p_{s+1}$ determines a
  partial flag -- a line and a plane -- passing through $p_0$. Choose coordinates $(x,y, z_1, z_2, \ldots , z_{n-2})$ in
  $\C^n$ centered at  $p_0$
 so that the line ($p_1$) is the x-axis
  and the plane is $x,y$ -- plane.   For convenience, set ${\bf z} = (z_1, z_2, \ldots z_{n-2}) \in \C^{n-2}$.
Consider  the locus of all
  points $p_{s+2}$ of type $R^sVT$ lying over $p_{s+1}$ and the corresponding curve
  germ set consisting of   $\gamma \in Germ(p_{s+2})$
  as $p_{s+2}$ varies over this locus.
      According to the KR computations   above (the case $R^s VT$ within the Proof of the codimension 2 theorem, theorem   \ref{thm:cod2class}) any such curve $\gamma$
  takes  the form
  \begin{equation}
  \label{eq:formgamma}
  \gamma = (t^3, a_1 t^{3s +1}  + a_2 t^{3s + 2} ,  {\bf b}_2 t^{3s +2} )  + O(t^{3s +3}),  a_1 \ne 0
  \end{equation}
  after a reparameterization.  In this expansion ${\bf b}_2$ is a vector in $\C^{n-2}$. It
  plays the role of the scalar $b_5$ in the previous proof.
  The $s$-fold prolongation of $\gamma$ is immersed, with fiber coordinates of  the form
  $$ u = a_1 n_1 t + a_2 n_2  t^2 ,  {\bf v} = {\bf b}_2 n_2  t ^2  \in \C^{n-2}.$$
  where $n_j = n_j (s) = (3s +j)(3(s-1) + j) \ldots (3 + j) j/ 3^s$.
  By the lemma, the 2-jet of $\gamma^s$ up to reparameterization uniquely determines the point $p_{s+2} = \gamma^{s+2} (0)$.
  The modifier `up to reparameterization' requires care.  Instead, we directly compute:
   the KR coordinates  one step over $p_{s+1}$ are   $d ^2 x /du^2, d^2 {\bf v} /du^2$
  and are given by $(0, {\bf b}_2)$ for our curve.    We see that the point $p_{s+2}$ is spatial
  if and only if ${\bf b}_2 \ne 0$, which defines an open dense set within $R^s VT$.
  All points in this set are realized by a curve equivalent to $(t^3, t^{3s+1}, t^{3s+3} , 0, \ldots , 0)$.

\medskip

\paragraph{\bf Case $R^s VVR$.}
Fix a point $p_{s+2}$ representing $R^s VV$.
 All such points are equivalent, being represented by a germ RL equivalent to $(t^3, t^{3s+2})$.
As per proposition \ref{partialflag}, the one-step projection  $p_{s+1}$ of $p_{s+2}$ determines a
  partial flag -- a line and a plane -- passing through $p_0$. Choose coordinates $(x,y, {\bf z})$ in
  $\C^n$ centered at  $p_0$
 so that the line ($p_1$) is the $x$-axis
  and the plane is $x,y$ -- plane.
Now, consider  the locus of all
  points $p_{s+3}$ of type $R^sVVR$ lying over $p_{s+2}$ and the corresponding
  curve germ set   of all   curves $\gamma \in Germ(p_{s+3})$
  as $p_{s+3}$ varies over this locus.
      According to the KR computations  (the case $R^s VV$ within the Proof of the codimension 2 theorem, theorem   \ref{thm:cod2class})) any such curve $\gamma$
  takes  the form
  $$\gamma = (t^3, a_2 t^{3s +2}  + a_3 t^{3s + 3} + a_4 t^{3s +4}  ,  {\bf b}_3  t^{3s +3}   + {\bf b}_4 t^{3s +4} ) +  O(t^{3s +5}),  a_2 \ne 0 $$
  after a reparameterization.  Note, in this expansion the  ${\bf b}_i$ are vectors    in $\C^{n-2}$.
  The non-vanishing of the  vector
 ${\bf b}_4$ will tell us whether or not the point is purely spatial.

  The  $s$-fold prolongation of $\gamma$ has fiber  coordinates
    $$ u = a_2 n_2  t ^2 + O(t^4) ,  {\bf v} = {\bf b_3} n_3  t ^3 +  {\bf b_4} n_4 t^{4}   \in \C^{n-2}.$$
  where $n_j = n_j (s) = (3s +j)(3(s-1) + j) \ldots (3 + j) j/ 3^s$.
  The $(s+1)$-fold prolongation has KR coordinates $u_1 = dx/du, {\bf v}_1 = d {\bf v} /du$
  and so:
  $$u_1 = {3 \over {2 a_2 n_2}} t  + O(t^2), {\bf v}_1
  =  {3 \over {2 a_2 n_2}} {\bf b_3} n_3  t  +  {4 \over {2 a_2 n_2}} {\bf b_4} n_4 t^{2}  $$

  To simplify our computation, we will  first show that we can assume   ${\bf b}_3 = 0$.
 Indeed,   the value of ${\bf b}_3$
 determines the location of the point $p_{s+2}$ in the fiber over $p_{s+1}$, which as we have
mentioned  can be placed arbitrarily using a symmetry, all $R^s VV$ points being equivalent.
  To see this explicitly, note that KR fiber coordinates for the $(s+2)$-fold prolongation of $\gamma$ are
  $u_2 =d u /du_1$, ${\bf v}_2 = d{\bf v}_1/du_1$.  It follows that ${\bf v}_2 = c {\bf b}_3$
  where $c$ is a nonzero constant involving $a_2$ and the $n$'s.
  A diffeomorphism of the form $(x,y, {\bf z}) \mapsto (x,y, {\bf z} - a x^{s+1})$
  kills the term ${\bf b}_3$ in $\gamma$ and consequently shifts the ${\bf v}_2$ coordinate to zero.
  Thus,   fixing $p_{s+2}$ is tantamount to assuming ${\bf b}_3 =0$.

  We now have,  fiber coordinates along $\gamma^{s+1}$  at level $s+1$ of the form
   $$u_1 = {3 \over {2 a_2 n_2}} t  + O(t^2),  {\bf v}_1  = c {\bf b_4}  t^{2}  $$
   where, as above $c$ is a nonzero constant involving $a$'s  and   $n$'s.
    By lemma \ref{lem:eq}, the 2-jet of $\gamma^{s+1}$
    mod reparameterization uniquely determines the point $p_{s+3} = \gamma^{s+3} (0)$.
A reparameterization $t \mapsto \lambda t + O(t^2)$ has the effect on the ${\bf v}_1$ term of
$c {\bf b_4}  t^{2}  \mapsto c \lambda^2  {\bf b_4}  t^{2} $.  In particular this term is nonzero if and only if
the original curve $\gamma$ is equivalent to the curve germ $(t^3, t^{3s+2}, t^{3s+4})$.
Consequently, the set of $p_{s+3}$'s lying over $p_{s+2}$ fall into two types: those for which
 ${\bf b_4} \ne 0$ and thus are purely spatial, and those for which ${\bf b}_4 =0$
 and so their germs contain planar curves (equivalent to $(t^3, t^{3s+2}, 0)$).
 Because the correspondence (jets) $\to$ points is algebraic,  the locus of  these purely spatial points
 is open and dense.

\medskip

\paragraph{\bf The case $R^s V R^m V RR$}.  We will leave this tedious case to the reader. Q.E.D.

\subsubsection{Conclusion: Proof of the Theorem \ref{stableClasses} }

 We start with the case $R^s VT$.  We follow the initial
 set up and coordinates as with  the proof of Theorem \ref{thm:3Dclasses} in subsection \ref{RsVT}.
Thus we  fix the point $p_{s+1} \in R^s VT$ and adapted coordinates to $p_{s+1}$ as  in that proof.
Consider the set of all curve germs $\gamma$ passing through $p_{s+1}$
such that $\gamma \in p_{s+2 + q}$ for some point $p_{s+2 + q} \in R^s VT R^q$
lying over $p_{s+1}$.  As per the earlier proof, $\gamma^{s}$ is immersed, and $\gamma$ has the form
of equation \ref{eq:formgamma}.
By the lemma, knowing the    jet, $j^{2+q} \gamma^s$, up to reparameterization,  is equivalent
to knowing the point  $\gamma^{s+2+q} (0) = p_{s+2 + q}$.
Once $q \ge 1$, we   can fix the parameterization
by the insistence that $x = t^3$, since we can `see' $j^3 \gamma$ as part of
$j^{q +2} \gamma^s$.
The key observation
is simply that $\gamma^s$ includes the information of $\gamma$.  Once
we know $j^{6s +2} \gamma$ we know which of the finite number of singularity classes $E_{6s+2,p, i}$
(and its degenerations)  from Arnol'ds
list the curve lies in.   So, take any $q \ge 6s$.  Then a point in $p_{s+ 2 + q} \in R^s VT R^q$
determines the $6s +2$ jet of $\gamma \in Germ(p_{s+ 2 + q})$, and consequently
the particular singularity class.  Since every such jet of the given form represents precisely one
class, we have established the desired bijection.

The proof  for the other case, $R^s VV \setminus R^s VL$ is quite similar and omitted. Q.E.D.

 \section{Death  of the Jet identification number. Birth  of the Jet Set.}

In the book \cite{monz:monster}  the notion of `jet identification number''
was introduced and was a crucial tool in many of the classification
and normal form results there.  We take a moment to explain why this
notion fails in dimension 3 or more, and what might be salvaged out of it.

  Recall   the symbol $=_{rep}$ means ``equal up to reparameterization''.
   \begin{defn} We say that the point $p \in \Pp^k (n)$
 has jet number  $r$ if there is a $\beta \in  J^r ( \C, \C^n)$
 such that $Germ(p) = \{ \gamma \text{ a curve germ} :  j^r \gamma =_{rep} \beta \}$.
 We call $\beta$ (mod reparameterization) the jet determining $p$.
 \end{defn}

The purpose of the jet identification number is
 to effectively reduce the size of $Germ(p)$ to a point, namely the $r$-jet $\beta$
 (mod reparameterization) occuring in the definition.
  We showed that when $n =2$  every regular point (point whose code ends in R)
 has a jet identification number.  However, for $n=3$
 points do not have a jet identification number.

Let us see what goes wrong  with the jet identification number in
three dimensions by dividing the jet identification number
definition into two parts.
\begin{defn} [Jet identification number]
There is a unique  integer $r$ such that for all $\gamma  \in Germ(p)$
\begin{itemize}
\item[(a)] if $\tilde \gamma \in Germ(p)$ then $j^r (\gamma) =_{rep} j^r (\tilde \gamma)$
\item[(b)]  if $j^r (\gamma) =_{rep} j^r (\tilde \gamma)$ then $\tilde \gamma \in Germ(p)$.
\end{itemize}
\end{defn}

In dimension $n=2$ all points with code ending in $R$, and in particular, the points with
code $RVR$ have a jet identification number.
We claim  any  points of type  $RVR$ in dimension 3 have no jet identification number.

All such points are equivalent, so, we may as well
work with  the point $p_* = c^3 (0)$ where
$c$ is the standard cusp, $(t^2, t^3, 0)$.
 Consider the deformation of $c$, given by the family of
  curves $c_{\alpha, a,b} = (t^2, \alpha t^3 + a t^4, b t^4).$
 They all represent the class $RV$.    Let's look at the points
 $p_{\alpha, a,b} = c_{a,b}^3 (0)$ by computing  the corresponding fiber  KR coordinates

$$\text{level 1}:  u =d y/dx = {{ 3  \alpha} \over 2} t + 2  a t^2  ,  v = dz/dx = 2 b t^2;$$

$$\text{level 2}:  u_2 =d x/du = \frac{4t}{3 \alpha} - \frac{16 t^2 a}{2 \alpha^2} + O(3) ,  v_2 = dv/du =  \frac{8 b t}{3 \alpha} -\frac{64 a b t^2}{9 \alpha^2}  + O(3) ;$$

$$\text{level 3}:  u_3 = \frac{8}{9 \alpha^2} - \frac{256 a t}{27 \alpha^3} + O(2),   v_3 = \frac{16 b}{9 \alpha^2} - \frac{128 a b t}{9 \alpha^3} + O(2).$$ \\

The point $p_*$ has $u_3, v_3$ coordinates $8/9, 0$.
Fixing the value of the $u_3$ coordinate to be $8/9$ fixes $\alpha = \pm 1$,
and the two values are related by reparameterization. A bit of thought  now  shows that
fixing $p_*$ fixes the 3-jet of $c$.  Thus  curves in $Germ(p_*)$ have the same 3-jet up
to reparameterization, so by (a) of the definition  we must have  jet identification number $r \ge 3$.
But the curves $c_{0,b}$ have the same 3-jet as $c_{0,0}$ and have
$p_{0,b} \ne p_{0,0}$  and so are not in $Germ(p_*)$.  By (b) of the definition, the jet identification
of $p_*$ cannot be $3$.
On the other hand, since   $p_{a,0} = p_{0,0}$   there are
curves with different  4-jet from  $c_{0,0}$ but still lying in  $Germ(p_*)$,
which shows by (a) that the jet-identification of $p_*$ cannot be $4$.  The jet identification cannot
be greater than 4 since item (b) of the definition holds for $r = 4$ and any $r \ge 4$.

We  conjecture that with the exception of Cartan points, there are no points in $\Pp^k (n)$,
$k > 1, n >2$ with a jet identification number.

To try to rescue something from the jet identification number in dimension 3 and higher  we   observe that
points in
$Germ(p_*)$ in the above example are characterized by $\alpha = \pm 1, b =0$
but $a$ arbitrary:  the three-jet is determined, modulo reparameterization
and {\it part} of the 4-jet.   Instead of a jet identification number, we get a {\it jet set},
together with  a relevant {\it interval} of jet numbers, here $3$ and $4$.

\begin{defn} [Jet interval] Fix a point  $p$ of the monster.
Suppose that there are integers $r, R$ such that for all $\gamma  \in Germ(p)$
\begin{itemize}
\item[(a)] if $\tilde \gamma \in Germ(p)$ then $j^{r} (\gamma) =_{rep} j^{r} (\tilde \gamma)$
\item[(b)]  if $j^R (\gamma) =_{rep} j^R (\tilde \gamma)$ then $\tilde \gamma \in Germ(p)$.
The maximum of the integers  $r$ will be denoted by $r_1$.
The minimum of the integers $R$   will be denoted by $r_2$.
If $r_1 \le r_2$ then we call $[r_1, r_2]$ the jet interval.
\end{itemize}
\end{defn}

We note in  that if an integer $r$ exists as per item (a), then any integer less than $r$
also works in (a).  And that if an integer $R$ works as in item (b) of the definition,  then any integer greater than $R$
works in (b).  We  believe that $r_1 \le r_2$ is a consequence of the definitions.

In the case of the cusp above,  $[r_1, r_2] = [3,4]$.

\section{Intersection Combinatorics of critical planes.}

We say that a collection $\Lambda_i , i \in I$ of linear hyperplanes
in a vector space is ``in general position'' if, for every
subset $J \subset I$ of indices whose cardinality is
less than or equal to the dimension of the vector space we have that  the
codimension of the subspace  $\bigcap _{i \in J} \Lambda_i$
is equal to the cardinality $|J|$ of $J$.

\begin{thm}  \label{criticalplanes}
 The critical planes through $p \in \Pp ^k (n)$ are in general position within
the vector space $\Delta_k (p)$
and there are at most $n$ of them.
\end{thm}

Recall the following language since it will  be central to  the proof.
The baby Monster originating from $p_k \in \Pp ^k (n)$
 is the tower of submanifolds  with induced distributions
 $(\Pp^j (F_k (p)), \delta_k ^j ) \subset (\Pp^{k +j} (n), \Delta_{k+j})$  obtained by prolonging
 the fiber   $F_k = \pi_{k, k-1}^{-1} (p_{k-1}) \subset \Pp ^k (n)$ through $p_k$ at level $k$.

The proof of theorem \ref{criticalplanes}, by induction on the ambient dimension $n$,
will rely  on
\begin{prop}  \label{babycritplanes} Within the baby Monster, $\Pp^j (F_k (p))$ the critical planes
are all of the form $\delta_k ^j \cap \delta_r ^s$ with
$r > k$ and $r+s = k+j$.
\end{prop}

 \paragraph{\bf Proof of the proposition.}
 For $c > 0$ the fiber of     $ \Pp^c (F_k (p))  \to  \Pp^{c-1} (F_k (p))$
 is $F_{k+c} \cap  \Pp^c (F_k (p))$. It follows that the baby monsters
 {\it within} the baby Monster originating from $p_k$ are obtained
 by intersecting $p_k$'s baby Monster with those originating at higher levels.
 Taking tangents yields the proposition.

We will  also be relying  on the following linear algebra lemma,
stated without proof.
 \begin{lem}\label{linalg}
 If $\Lambda_0, \Lambda_1, \ldots , \Lambda_s $ is  a collection of  linear hyperplanes
 in $V$ such that collection $\Lambda_0 \cap \Lambda_1,  \ldots\Lambda_0 \cap \Lambda_s$
 is  in general position within $\Lambda_0$ then  the original collection
   is in general position within $V$
 \end{lem}

\paragraph{\bf Proof of theorem \ref{criticalplanes}.}

 By induction on $n$.   The base of the induction is  the case $n=2$ of the theorem,
 for  curves in the plane
and  was proved in the book \cite{monz:monster} where it was central to the development.

The inductive hypothesis is
\noindent {\it   hypothesis} $H_n$:  The critical planes
 through a point of $\Pp^k (n)$ are in general position and are
 less than or equal to $n$ in number.
 \smallskip
 \noindent
 Assuming $H_{n-1}$, we will prove $H_n$.

 To this end let $p \in \Pp^k (n)$.
 Let $k_0 \le k$ be the smallest integer
 such that a critical hyperplane $\delta_{k_0} ^{j_0}$ originating from
 level $k_0$ passes through $p$.  (Necessarily $j_0 + k_0 = j + k$.)
 Then the  successive levels  of   the baby Monster   arising from   level $k_0$
 pass
 through successive points $p_{k_0 +i}$, $i = 0, 1, \ldots , k-k_0 $  in the tower under $p$.
 By   proposition \ref{babycritplanes}, the critical planes of this baby Monster are of the
 form $\delta_{k_0} ^{i} \cap \delta_r ^s$, $r > k_0$.  By the inductive hypothesis
 $H_{n-1}$
applied to the $(n-1)$-manifold $F_k$ we have that these planes are in general position
and there are no more than $n-1$ of them.   By the linear algebra lemma \ref{linalg}
the collection $\delta_{k_0} ^{i}$ together with  $\delta_r ^s$ is in general position and there
are no more than $n$ of them. Q.E.D.

\medskip

 Theorem \ref{criticalplanes} asserts that there are at most $n$ critical planes passing through any
 point.  But we can often do better than this.  Through a Cartan point there   passes  exactly
 one  critical hyperplane, namely the vertical hyperplane.  Through an $A_{2k}$ point , i.e. one with code
 $R^k V$ there pass two critical planes, the vertical hyperplane, and the critical hyperplane arising
 from the baby monster one level down.  To get a bound in the general situation
let $p \in \Pp^k (n)$ let $\omega = \omega_1 \ldots \omega_k$
be its RC code.

\begin{defn}  The length of the critical tail for $p$ is equal to
zero if $\omega_k = R$. If $\omega_k = C$ then this length is
 the number of consecutive letters ending with
 $\omega_k$  which are C's
 \end{defn}

 \begin{prop} The number of critical planes through $p$ is less than or equal to
 $1$ plus the  length of the critical tail through $p$.
 \end{prop}

\paragraph{\bf Proof.} Follow  how the critical planes through $p$
 can arise out of baby monsters.  Observe that if a projection $p_i$ of $p$
 is a regular point, (that is, the corresponding line one level down is
 a regular direction)  then no critical planes can arise out of the  baby monster
 one level down, or further, since these baby monsters cannot pass through $p_i$, it being a regular point. Q.E.D.

\subsection{RVL is non-planar}

We give two applications of
the incidence relations for critical planes.
Recall that  we proved that any curve representing
the   class $R^s VL R^q$ is of the form  $(t^4, t^{4s+1}, t^{4s+2}, \ldots )
+ O(t^{4s+4})$.  This jet is that of a non-planar curve:  there is no smooth
surface which contains the curve near $t=0$.  It follows that no   point in this RVT class can be touched
by a planar curve.   We    give an alternative   `synthetic' proof of the  non-planarity
of these points,  a proof   which   simultaneously   establishes the non-planarity
of points in  any RVT class having the letter `L'.
 We then   generalize this
theorem to higher-dimensional curves.

\begin{thm}
Let $\zeta = (p,\ell)$ be a point  of the monster at level $k+1$
for which the  line   $\ell$
lies in the  intersection of   two     critical planes through the point $p$ at level $k$.
Then $\zeta$ cannot be touched  by the prolongation of a  planar curve.
\end{thm}

\paragraph{\bf Proof.}  It suffices to prove that
if $\gamma$ is planar then its prolongations $\gamma^k$
are never tangent to a line lying in two critical planes. Since
the condition of planarity and the condition of being a line contained in two
critical planes are diffeomorphism invariant conditions, we may use   symmetry considerations,
and take the   surface containing
$\gamma$ to be   the xy-plane   $F_0$ sitting inside $\R^3$.

We can now work   entirely in the $n=3$ context.  The trick is to treat  $F_0$ in  the same way
as we treated  fibers  in the baby monster construction,  thinking of it
as a ``level 0 fiber'', and prolonging it so as
to obtain the ``level 0 baby monsters'' , the  submanifolds-with-distribution $(\Pp ^k (F_0), \delta_0^k)$
within $\Pp ^k (3)$. We have that   $\delta_0^k\subset \Delta_k$ is a 2-plane to which $\gamma^k$
   is tangent.
Throw  the $\delta_0^k$ in to the collection of critical hyperplanes and observe that the  resulting larger collection of
  critical hyperplanes
is still in general position, as is seen by going back over the proofs of
Theorem \ref{criticalplanes} and   Proposition \ref{babycritplanes}, viewing $F_0$ as another fiber.   Now,
a given line $\ell$ in a 3-space can be contained in at most 2 two-planes out of a collection of
two-planes which are in general position.  Our line $\ell$, the tangent line  to $\gamma^k$,
  lies in   $\delta_0 ^k$ since $\gamma \subset F_0$.
It follows that  $\ell$ can lie in at most one other critical hyperplane.
Q.E.D. \\

Recall that a curve  has embedding dimension  $d$ if it lies
in a smooth $d$-dimensional submanifold.

A nearly identical proof to that just given yields,
 \begin{thm}
Let $\zeta = (p,\ell)$ be a point  of the monster
for which the  line   $\ell$
lies in the  intersection of  $d-1$ critical hyperplanes through the point $p$.
Then $\zeta$ cannot be touched  by the prolongation of a curve with embedding dimension less than $d$.
\end{thm}

This is a special instance of E. Casas-Alvero's projection theorem in the language of the prolongation tower. See (\cite{cas:proj}, especially  pp. 318).

\subsection{Orbit counts in the spatial case.}
\begin{figure}
  \includegraphics[width=7.0in]{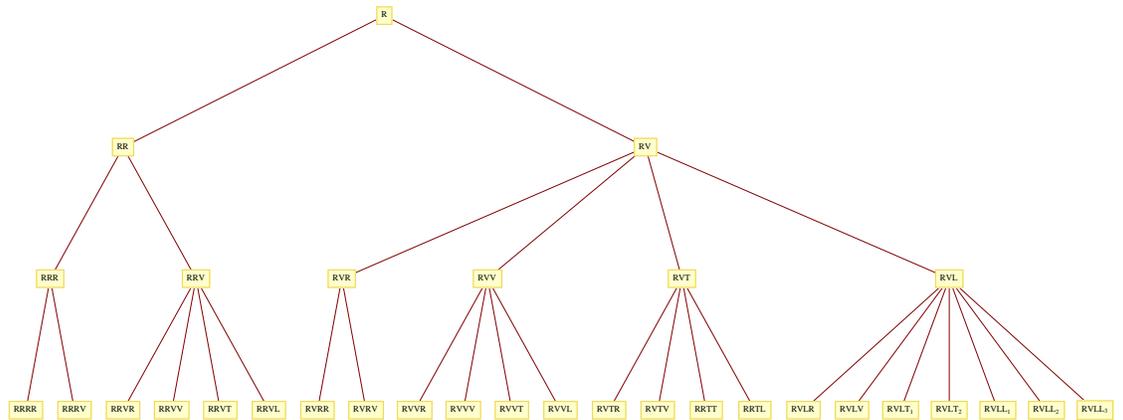}\\
  \caption{Classification of geometric orbits in $P^k(3)$ for $k\leq 4$.}\label{fig:orbit-branching}
\end{figure}

The introduction of the class `L' allows us to further split  the RVT classes into many more classes.
Theorem \ref{criticalplanes} asserts that at most   three critical planes
can pass through a point when   $n=3$. The maximum of $3$  is realized  if and only if the point is of
`L' type.  Then its three critical planes are  in general position, and each pair gives rise to another `L' direction, and hence a new type of L point one level up.  Therefore we have 3 directions of type `L', and two distinct types of tangency directions since their corresponding baby monsters are born in different levels.   An immediate conclusion of this simple reasoning is that there are at least 7 geometrically distinct types of directions passing through a `L' point : regular, vertical, 2 types of tangency, and 3 types of `L' directions. Adding together these contributions level-by-level we compute a {\it lower bound} for the number of orbits in the spatial case $n=3$ summarized in table \ref{table:orbit-counting}. The tree graph in figure \ref{fig:orbit-branching} summarizes this `branching' of geometric classes in first four levels of the extended monster tower.

\begin{table}\label{table:orbit-counting}
\begin{tabular}{|c|c|c|}
  \hline
  level & planar orbits(sharp) & spatial orbits \\ \hline
  first & 1 & 1 \\ \hline
  second & 2 & 2 \\ \hline
  third & 5 & 6 \\ \hline
  fourth & 13 & 23 \\ \hline
  fifth & 34 & 98 \\
  \hline
\end{tabular}
\caption{Orbit counting comparison between planar and spatial cases. The planar number are exact. In the spatial case the numbers are lower bounds.}
\end{table}

\medskip

\section{Semple tower = Monster Tower}

Algebraic geometers, beginning with Semple (\cite{semp:tower}) in the 1950s, have
been working on the Monster Tower, which they now call the Semple Tower.
See (\cite{ken:cnct}) and (\cite{LJ:semple}).
Semple's original tower concerned the planar case, and has base
$\bP^2$ rather than our $\C^2$.

M. Lejeune-Jalabert  makes a
particularly  beautiful use of the Tower  in (\cite{LJ:semple}) in order to generalize
F. F. Enriques (\cite{enr:hist})  famous formula relating
the multiplicities of points on consecutive blow-ups of singular plane curves.
  M. Lejeune-Jalabert
 generalized this formula to the case of  curves in $\C^n$.
We tested our development of the Monster
 against her results to  obtain an alternative  derivation of her Enriques' formula.

\begin{thm}[Enriques]
Given a germ $\gamma:(0,\C)\rightarrow (0,\C^n)$, let $S = \{\gamma^{i}(0)\}_{i=0}^{\infty} = \{p_i\}$. Then,
\begin{eqnarray*}
\mbox{multiplicity of }p_i = \sum_{p_j \small{\mbox{ proximate to }} p_i} \mbox{ multiplicity of } p_j
\end{eqnarray*}
\end{thm}

The statement of the theorem uses  the notion of   points at different levels being {\em proximate}.
{\em Proximity} is a classical notion in the algebraic theory of curves, applied to sequences
of `infinitely near' points obtained from classical blow-ups of the curve.  Semple (\cite{semp:curv} or \cite{lj:semple}),
and then Lejeune altered the notion  so as to fit  the Nash blow-up (= prolongation) of the
curve as it fits within the Semple tower.   The
following definition  is equivalent to Lejeune's  (op. cit. \cite{lj:semple}) .\\
\begin{defn}A point $p\in \mathcal{P}^k(n)$ is said to be {\em proximate} to a point $q$, $q = \pi_{k,j}(p)$ under $p$ in the tower
if  either:
\begin{enumerate}
\item $j = k-1$
\item there a vertical curve $\sigma$ through $q$ whose prolongation (sufficiently many times)
passes through $p$.
 \end{enumerate}
 \end{defn}
From our structure theorem on the critical hyperplanes through a given point, any point in $p\in \mathcal{P}^k(n)$ is proximate to at most $n$ points. (Sitting at lower levels. )  \\

By keeping track of the multiplicity of consecutive prolongations of a well-parameterized singular germ $\gamma$
one can show that after sufficiently many prolongations the prolonged curve is a {\em regular germ}.  This fact yields an alternative proof of the following theorem of Nobile \cite{nob:nash}:
\begin{thm} [Nobile; Castro] \label{Nobile}
 Every well-parameterized curve germ has a  finite regularization level.
\end{thm}
We shall present our demonstration in a different note.

\medskip

\section{Open Problems.}

We end with some open problems.
Throughout    $p \in \Pkn$
and $\gamma$ is a  non-constant well-parameterized curve germ in $\C^n$.

\subsubsection{ On simplicity}
\paragraph{Q1. :} Is $p$ is tower simple if and only if every $\gamma \in Germ(p)$ is simple?
\medskip
\paragraph{Q2. :} Is  $\gamma$ is simple if and only if all the points  $\gamma^j (0)$   are simple?
\medskip
There are also stable versions of these two questions.  To formulate them, we
use  the embeddings $\C^n \to \C^{n+1}$ to obtain  embeddings $\Pkn \to \Pp^{k} (n+1)$ which
take distribution into distribution.     The following diagram (``the Russian doll") may be helpful:
\def\mapright#1
{
\smash{
\mathop{\longrightarrow}\limits^{#1}}
}
\def\mapdown#1
{
\downarrow{\rlap{#1}}
}

$$
\begin{matrix}
\mapdown{} &   & \mapdown{}  &  & \ldots & & \mapdown{} & \ldots  \cr
\Pp ^k (2)   & \mapright{} & \Pp^k (3)  &  \mapright{} & \ldots & \mapright{} &\Pp^k (n) & \mapright{}  \cr
\mapdown{} &   & \mapdown{}  &  & \ldots & & \mapdown{} & \ldots  \cr
\vdots &   & \vdots   &  & \vdots & & \vdots & \ldots  \cr
\mapdown{} &   & \mapdown{}  &  & \ldots & & \mapdown{} & \ldots  \cr
\Pp ^1 (2)   & \mapright{} & \Pp^1 (3)  &  \mapright{} & \ldots & \mapright{} &\Pp^1 (n) & \mapright{}  \cr
\mapdown{} &   & \mapdown{}  &  & \ldots & & \mapdown{} & \ldots  \cr
\mathbb A^2  &\mapright{} & \mathbb A^3 & \mapright{} & \ldots  & \mapright{} & \mathbb A^n & \mapright{}  \cr
\end{matrix}
$$

 Arnol'd used the standard embeddings $\C^n \to \C^{n+1}$ to define the notion of
a ``stably simple curve singularity'':  one which is simple independent of the embedding dimension.
He listed all these stably simple curves.

\medskip

\paragraph{Q3:} Does iterated prolongation induce a  bijection between
Arnol'ds stably simple curves and the R-stabilizations of stably simple points?

\subsubsection{On discrete invariants attached to points}

\paragraph{{\bf Discrete invariants.}}
\medskip
\paragraph{Q4. :}  For $p \in \Pkn$ is it true that $mult(\gamma)$
is independent of the choice of $\gamma \in Germ(p)$?
\medskip

If so, we would call this number the `multiplicity of $p$'.
(If not we could take the minimum of the multiplicities, but that
would be less satisfactory.)

Instead of the multiplicity of a curve,
we could take any discrete invariant $\Lambda(\gamma)$ of  curve germs
and ask, for  given $p$ is the value $\Lambda (\gamma)$
constant for all $\gamma \in Germ(p)$?  If `yes' we would then say    that
the invariant $\Lambda$ is well-defined for $p$.   Possible invariants are
the semi-group and the parameterization number.    Again, if the invariant
is not well-defined, but sits within a partially ordered set, perhaps we can take a
minimum of its values on $\gamma  \in Germ(p)$ to  get an invariant
of the point.   In addressing these questions, and the ones around them,
some  substitute for the jet identification number, so useful in dimension 2, but
debunked in dimension 3 and greater, will be of great help.

\medskip
\paragraph{Q5. :}
Does every regular point (one whose RVT code ends in R) in $\Pkn$,  $n > 2$
have a jet interval, with associated jet set, as defined above?
\medskip

\subsubsection{ Curve-to-Point philosophy.  Continuity of the prolongation-evaluation map.}

 Write $Germ_w ^{\ast} \subset Germ^{\ast}$ denote  the space of all {\it well-parameterized}
curve germs $(\C, 0) \to (\C^n, 0)$.  We  have explicitly excluded the constant curve.
Write $\Pp_0$ for the fiber of the infinite Monster over $0 \in \C^n$, this being the direct
limit as $k \to \infty$ of the $\Pkn$.
A point of $\Pp_0$ is an infinite sequence $(0, p_1, p_2, \ldots, p_i, p_{i+1}, \ldots)$
with $p_k \in \Pkn$ and $\pi_{k, i} (p_k) = p_i$.
We have the prolongation map
$$Prol:  Germ_w ^{\ast} \to \Pp_0$$ by sending
$\gamma$ to $Prol(\gamma) = (\gamma(0), \gamma^1 (0), \ldots, , \gamma^k (0), \ldots)$.
By theorem \ref{Nobile}, $\gamma^i$  is regular for $i$ sufficiently large,
so that all the  the $\gamma^i (0)$ are regular points.
It follows that the range of $\Pp_0$  is contained in the subset
$\Pp_0 ^R$  of points which are   eventually
regular.

A coordinate computation with immersed curve germs shows that the image of $Prol$
is not all of $\Pp_0 ^R$: some kind of ``growth'' conditions on KR coordinates
are also required if the domain of $Prol$ consists of analytic functions.
For consider the     coefficients $b_j$ of the coordinate functions $x_a $
along the curve.  Being analytic functions, these coefficients
satisfy  a bound $|b_j|  \le C^j$ for some constant $C$ depending only
on $\gamma$.   Preparing the coordinates and parameterization
so   $x(t) = t$, the    KR coordinates  of $Prol(\gamma)$ at
level $j$, written $(u_j, v_j, \ldots)$, have the form  are of the form $d^j x_a /dx^j$, etc, and so satisfy the bounds
$|u_j| \le j! C^j$.

\medskip

\paragraph{Q6. :} What is the image of $Prol$?

\medskip

Let us denote  this image by $\Pp_0 ^{an} \subset \Pp_0 ^R$,
the `an' being for analytic.

It  seems quite  impossible  to make  $Prol$ into a continuous map relative to   any reasonable
topology on the space of curve germs.  To see the problem,  take $N$ large and consider the curve
$\gamma = (x(t), y(t) = (t^N, t^{N+1})$.  Then $\gamma^1 (0) = (0,0,0)$ in standard coordinates
where the last coordinate represents $dy/dx$.  But, for any $A \in \C$, and $r < N$
we can find an arbitrarily $C^r$-small perturbation $\tilde \gamma$ of $\gamma$
with $\tilde \gamma^1 (0) = (0,0,A)$.
However, restricted to curves of multiplicity $1$(immersed curves) $Prol$ is beautifully
continuous.

\medskip

\paragraph{Q7. :}  Is there some discrete curve  invariant $\Lambda$ such as multiplicity,
such that    $Prol$ is continuous when   restricted to the class of all  curves  of constant $\Lambda$?

As a kind of converse  to the previous question we ask .
\medskip
\paragraph{Q8. :}  If $\gamma$ is a well-parameterized curve germ  with $\gamma^j$ regular,
then, is there an $r \ge j$ such that a \nbhd base for $\gamma^j (0) \in \Pp ^j (n)$
induces a $C^r$ \nbhd base  the curve germ $\gamma$?
An affirmative answer would yield a potentially powerful technical and conceptual tool.
\medskip

 \bibliography{1nGb}
 \bibliographystyle{amsplain}




\end{document}